\documentclass[psamsfonts]{article}

    \usepackage{latexsym}
    \usepackage{amsfonts}
    \usepackage{amsmath}
    \usepackage{amssymb}
    \usepackage{amscd}
    \usepackage{amstext}
    \usepackage{amsthm}

 \title{Minimal Homeomorphisms and Approximate Conjugacy in Measure
 \thanks{Supported by National Science Foundation
of USA. \,\,\,\,\,\,\,\,
 AMS 2000 Subject Classification Number: Primary 37A55 and 46L35.\,\,\,\,\,\,\,\,\,\,\,
Key Words: Simple $C^*$-algebras, Approximate Conjugacy\protect\\}}

 \author{
 Huaxin Lin\\
    Department of Mathematics\\
    University of Oregon\\
    Eugene, Oregon 97403-1222\\
    }
 \date{}
  
\begin{document}
\maketitle

    \newcommand{\CA}{$C^*$-algebra}
    \newcommand{\SCA}{$C^*$-subalgebra}
\newcommand{\aue}{approximate unitary equivalence}
    \newcommand{\ayue}{approximately unitarily equivalent}
    \newcommand{\mops}{mutually orthogonal projections}
    \newcommand{\hm}{homomorphism}
    \newcommand{\pisca}{purely infinite simple \CA}
    \newcommand{\andeqn}{\,\,\,\,\,\, {\rm and} \,\,\,\,\,\,}
    \newcommand{\QED}{\rule{1.5mm}{3mm}}
    \newcommand{\morp}{contractive completely
    positive linear map}
    \newcommand{\asmorp}{asymptotic morphism}
    \newcommand{\arrow}{\rightarrow}
    \newcommand{\tdsum}{\widetilde{\oplus}}
    \newcommand{\pa}{\|}  
    \newcommand{\ep}{\varepsilon}
    \newcommand{\id}{{\rm id}}
    \newcommand{\aueeps}[1]{\stackrel{#1}{\sim}}
    \newcommand{\aeps}[1]{\stackrel{#1}{\approx}}
    \newcommand{\dt}{\delta}
    \newcommand{\yu}{\fang}
    \newcommand{\ca}{{\cal C}_1}
\newcommand{\Ad}{{\rm ad}}
    \newcommand{\tr}{{\rm TR}}
    \newcommand{\N}{{\bf N}}
    \newcommand{\C}{{\bf C}}
    \newcommand{\Aut}{{\rm Aut}}
    \newcommand{\Tand}{\,\,\,\text{and}\,\,\,}
\newcommand{\T}{{\mathbb T}}
\newcommand{\R}{{\mathbb R}}
    \newtheorem{thm}{Theorem}[section]
    \newtheorem{Lem}[thm]{Lemma}
    \newtheorem{Prop}[thm]{Proposition}
    \newtheorem{Def}[thm]{Definition}
    \newtheorem{Cor}[thm]{Corollary}
    \newtheorem{Ex}[thm]{Example}
    \newtheorem{Pro}[thm]{Problem}
    \newtheorem{Remark}[thm]{Remark}
    \newtheorem{NN}[thm]{}
    \renewcommand{\theequation}{e\,\arabic{section}.\arabic{equation}}
    \newcommand{\rforal}{{\rm\,\,\,for\,\,\,all\,\,\,}}
\newcommand{\Z}{{\mathbb Z}}
\newcommand{\Q}{{\mathbb Q}}
\newcommand{\di}{{\rm dist}}
       \newcommand{\af}{\alpha}
       \newcommand{\bt}{\beta}
       \newcommand{\Om}{\Omega}
       \newcommand{\gm}{\gamma}
      \newcommand{\sm}{\sigma}
      \newcommand{\ud}{\underline}
      \newcommand{\beq}{\begin{eqnarray}}
     \newcommand{\eneq}{\end{eqnarray}}

\begin{abstract}
Let $X$ be an infinite  compact metric space with finite covering dimension. Let $\af,\bt: X\to X$ be
two minimal homeomorphisms.        Suppose that the range of $K_0$-groups of both crossed products are
dense in the space of real affine continuous functions.
 We show that
$\af$ and $\bt$ are approximately conjugate uniformly in measure if and only if they have affine
homeomorphic  invariant probability measure spaces.
\end{abstract}

\section{ Introduction}

\hspace{0.2in} Let $X$ be a compact metric space and let $\af, \bt: X\to X$ be two minimal
homeomorphisms. If $X$ has infinitely many points, then the associated crossed product \CA s
$C(X)\rtimes_{\af}\Z$ and $C(X)\rtimes_{\bt}\Z$ are unital
separable  simple \CA s. It was proved by J. Tomiyama (\cite{T}) that $\af$ and $\bt$ are conjugate if
there is a *-isomorphism
$\phi$ from $C(X)\rtimes_{\af}\Z$ onto $C(X)\rtimes_{\bt}\Z$ which
maps $C(X)$ onto $C(X).$ On the other hand, T. Giordano, I. Putnam and C.  Skau (\cite{GPS}) showed
that two minimal Cantor systems are topological orbit equivalent if and only if the tracial range of
$K_0(C(X)\rtimes_{\af}\Z)$ is unital order isomorphic to that of $K_0(C(X)\rtimes_{\bt}\Z).$ Both
results show the strong connection between \CA\, theory and minimal dynamical systems. There are also
many other studies on the interplay between \CA\, theory and minimal dynamical systems.  With the
development of classification of simple amenable \CA s (example, \cite{EG} and see also
\cite{Lnduke}), one hopes that there will be  further applications of \CA\, theory to the minimal
dynamical systems as well as the application of theory of topological dynamics to the \CA\, theory. We
believe that $K$-theory or tracial information of the crossed products may play more interesting roles
in the study of minimal dynamics than what we currently know.
 Conjugacy in minimal dynamical systems is
certainly a very strong relation. Orbit equivalence studied in
\cite{GPS}  is much weaker than conjugacy in general. However for connected spaces orbit equivalence
often implies the conjugacy. Approximate conjugacy in minimal dynamical systems have been introduced
and studied recently (see \cite{Lncd}, \cite{LM1}, \cite{M}, \cite{LM2} and \cite{LM3}) which are
closely related to the classification of amenable simple \CA s.  In this paper, we study a much weaker
equivalence relation, namely, uniform approximate conjugacy in measure.

Roughly speaking, two minimal homeomorphisms $\af$ and $\bt$ are approximately conjugate uniformly in
measure if there exists a sequence of Borel isomorphisms $\gm_n: X\to X$ such that
$\gm_n^{-1}\af\gm_n$ converges to $\bt$ and $\gm_n\bt\gm_n^{-1}$
converges to $\af$ in measure uniformly on the set of
$\bt$-invariant measures and the set of $\af$-invariant measures,
respectively. We also require that $\{\gm_n\}$ eventually preserves measures in a suitable sense.
Moreover, $\{\gm_n\}$ and $\{\gm_n^{-1}\}$ should be continuous on some (eventually dense) open
subsets of $X.$ The precise definition is given in \ref{Dapp1}.

One should not expect that $\{\gm_n\}$ converges in any reasonable sense in general. So it is
important that maps $\gm_n$ have certain consistency. Since homeomorphisms do not preserve measures,
it is crucial to know, among other things, that
$\mu(\gm_n(E))\not\to 0,$ or $\nu(\gm_n^{-1}(E))\not\to 0$ for any
open sets $E$ and $\gm_n(E),$ or $\bt$-invariant measure, respectively. Therefore it should be
regarded as a crucial part of the definition that $\gm_n$ as well as
$\gm_n^{-1}$ preserve measures in certain sense.

Let $A_\af=C(X)\rtimes_{\af}\Z$ and $A_\bt=C(X)\rtimes_{\bt}\Z.$ Denote by $T_\af$ and $T_\bt$ the
compact convex sets of
$\af$-invariant probability Borel measures and $\bt$-invariant
probability Borel measures, respectively.
 Under
the assumption that
$\rho(K_0(A_\af))$ is dense in
$Aff(T(A_\af))$ and
$\rho(K_0((A_\bt)))$ is dense in $Aff(T(A_\bt))$ (see \ref{P1}
(4), \ref{P4} and \ref{LP} below),  we prove that if $T_\af$ and
$T_\bt$ are affine homeomorphic, then $\af$ and $\bt$ are
approximately conjugate uniformly in measure (see Theorem \ref{MT1} below). If both $T_\af$ and
$T_\bt$ have only finitely many extremal points, this condition simply says that $T_\af$ and
$T_\bt$ have the same number of extremal points. So this condition
requires little. Therefore, one should not regard uniform approximate conjugacy in measure as a strong
relation. To the contrary, we would like to emphasis that
 two
minimal homeomorphisms could be easily approximately conjugate uniformly in measure. In particular, if
both $\af$ and $\bt$ are uniquely ergodic, then they are always approximately conjugate uniformly in
measure.
However, we would also like to point out that the cases that are most interesting here are
the cases that
$\af$ and $\bt$ have rich invariant measures.
Given an affine
homeomorphism $r: T_\af\to T_\bt,$ Theorem \ref{MT1} says that $r$ can be induced by a sequence of
Borel equivalences $\{\gm_n\}$ of
$X$ for which $\gm_n^{-1}\af\gm_n$ converges to $\bt$ and
$\gm_n\bt\gm_n^{-1}$ converges to $\af$ in measure uniformly ( not
just for each $\mu\in T_\af$ and $\nu\in T_\bt$).

In a subsequent paper, we will discuss a much stronger measure theoretical version of approximate
conjugacy which is closely related to the Giordano-Putnam-Skau orbit equivalence for Cantor minimal
systems.

The paper is organized as follows. Section 2 lists a number of notation and facts used in this paper.
Section 3 gives a  versions of uniform Rohlin property  for
dynamical systems with mean dimension zero. Section 4 contains a number of technical
lemmas which will be used in the proof of the main result of the paper. Section 5 discusses the notion
of uniform approximate conjugacy in measure and present the proof of the main result (Theorem
\ref{MT1}). Finally, section 6 gives a few concluding remarks.

\vspace{0.2in}

{\bf Acknowledgement} This work is partially supported by a grant from National Science Foundation of
U.S.A. This work was initiated when the author was visiting East China Normal University in the Summer
2004. It was partially supported by Shanghai Priority Academic Disciplines.

\section{Preliminaries}

\begin{NN}\label{P1}
\vspace{0.2in}
 \indent
 {\rm (1) If $k$ is a positive integer, $M_k$ is the full matrix
algebra over ${\mathbb C}.$ Denote by ${\rm Tr}$ the standard trace on $M_k$ and ${\rm tr}$ the
normalized trace on $M_k.$

\vspace{0.1in}

 (2) Let $A$ be  a \CA. Denote by $T(A)$ the tracial
state space of $A.$ If $A$ is stably finite, $T(A)\not=\emptyset.$ It $\tau\in T(A),$ we will also use
$\tau$ for $\tau\otimes Tr$ on
$M_k(A),$ $i=1,2,....$

\vspace{0.1in}

(3) Let $Aff(T(A))$ be the space of all real affine continuous functions on $T(A).$ Let $a\in
A_{s.a.}.$ Denote by $\hat{a}$ the real affine continuous function defined by $\hat{a}(\tau)=\tau(a)$
for $\tau\in T(A).$

\vspace{0.1in}

(4)  Denote by $\rho_A: K_0(A)\to Aff(T(A))$ the order \hm\, induced by $\hat{p}$ for projections
$p\in M_k(A),$ $k=1,2,....$
 We often use $\rho$ if the \CA\, $A$ is understood.
}
\end{NN}

\begin{NN}\label{P2}

{\rm (5) Let $X$ be a compact metric space. We say $X$ has finite dimension if $X$ has finite covering
dimension.

\vspace{0.1in}

(6) Let $h: C(X)\to A$ be a contractive positive linear map. Suppose that $t$ is a positive linear
functional of $A.$ Then
$t\circ h$ gives a positive linear functional of $C(X).$ We will
use $\mu_{t\circ h}$ for the positive Borel measure on $X$ induced by the positive linear functional
$t\circ h.$ }
\end{NN}

\begin{NN}\label{PM}
{\rm (7)  Let $X$ be a compact metric space and $\af: X\to X$ be a homeomorphism. Recall that $\af$ is
minimal if $\af$ has no proper
$\af$-invariant closed subset, or, equivalently, for each $x\in
X,$ $\{\af^n(x): n=0,1,2,...\}$ is dense in $X.$

\vspace{0.1in}

 (8) Let $X$ be a compact metric space and $x\in X.$ The point $x$ is said to be a
 {\it condensed point} if every open neighborhood of $x$ contains uncountably many
 points of $X.$

\vspace{0.1in}

 (9) If $X$ has infinitely many points and $\af$ is
minimal, then the cross product $C(X)\rtimes_\af \Z$ is a unital simple \CA. We will use $A_\af$ for
$C(X)\rtimes_{\af} \Z.$

   In this case, $X$ has no isolated points and every point of $X$ is condensed.

\vspace{0.1in}

 (10) Denote by $j_\af: C(X)\to A_\af$ the usual
embedding. Denote by $u_\af$ the implementing unitary in $A_\af$ such that
$$
u_{\af}^*j_\af(f)u_\af=j_\af(f\circ \af)\rforal
f\in C(X).
$$

\vspace{0.1in}

 (11) Denote by $T_\af$ the space of all
$\af$-invariant probability Borel measures on $X.$ If $\mu\in
T_\af,$ then it induces a tracial state $\tau_{\mu}$ so that
$$
\tau_{\mu}(j_\af(f))=\int f d\mu
$$
for all $f\in C(X).$ On the other hand, if $\tau\in T(A_\af),$ then $\mu_{\tau\circ j_\af}$ gives a
measure in $T_\af.$ This measure will be denoted by $\mu_{\tau}.$

The reader may notice that we  do not always attempt to distinguish the convex sets $T_\af$ from
$T(A_\af).$

 }
 \end{NN}

\begin{NN}\label{P3}

{\rm  (12) Let $A$ and $B$ be two \CA s. By a \hm\, $h: A\to B,$ we mean a $*$-\hm\, from \CA\, $A$ to
$B.$ Suppose that both $A$ and $B$ are unital and stably finite. We say that $r: Aff(T(A))\to
AffT(B))$ is a unital order \hm\, if $r$ is an order \hm\, and $r(\hat{1_A})=\hat{1_B}.$ The \hm\, $r$
is said to be an order isomorphism if $r$ is a bijection and $r^{-1}$ is an also order \hm.

Suppose that an affine continuous map $r: Aff(T(A))\to Aff(T(A))$ is a unital order isomorphism.
Denote by $r_{\natural}: T(B)\to T(A)$ the affine continuous map induced by
$r_{\natural}(\tau)(a)=r(\hat{a})(\tau)$ for all $a\in A_{a.s}$
and $\tau\in T(B).$ If $r$ is a unital order isomorphism, then
$r_{\natural}$ is an affine homeomorphism.

On the other hand, if $\lambda: T(A_\bt) \to T(A_\af)$ is an affine homeomorphism, then one obtains a
unital order isomorphism
$\lambda^{\sharp}: Aff(T(A_\af))\to Aff(T(A_\bt))$ by
$\lambda^{\sharp}(a)(\tau)=a(\lambda(\tau))$ for all $a\in
Aff(T(A_\af))$ and $\tau\in T(A_\af).$

\vspace{0.1in}

 (13) If $\phi: A\to B$ is a \hm\, we will use
$\phi_*: K_*(A)\to K_*(B)$ for the induced map on $K$-theory.

\vspace{0.1in}

 (14) Let $A$ and $B$ be two \CA s and $\phi: A\to
B$ be a \morp. Suppose that ${\cal G}$ is a subset of $A$ and
$\dt>0.$ We say $\phi$ is ${\cal G}$-$\dt$-multiplicative if
$$
\|\phi(ab)-\phi(a)\phi(b)\|<\dt
\rforal  a,\, b\in {\cal G}.
$$

\vspace{0.1in}

 (15) Let $\phi: C(X)\to A$ be a \hm. We say that
$\phi$ has finite dimensional range if the image of $\phi$ is
contained in a finite dimensional \SCA\, of $A.$ If $\phi$ has finite dimensional range, then there
are finitely many points
$\{x_1,x_2,...,x_m\}\subset X$ and pairwise orthogonal projections
$p_1,p_2,...,p_m$ in $A$ such that
$$
\phi(f)=\sum_{i=1}^m f(x_i)p_i
\rforal  f\in C(X).
$$

\vspace{0.1in}

 (16) Let $A$ be a unital simple \CA. We write
$TR(A)=0$ if $A$ has tracial rank zero. For the definition of
tracial rank zero, we refer to \cite{Lntr0} or 3.6.2 of \cite{Lnb}. A unital simple \CA\, with tracial
rank zero has real rank zero, stable rank one and weakly unperforated $K_0(A)$ (see \cite{Lntr0}).

}
\end{NN}

\begin{NN}\label{P4}

{\rm (17)}\,
 {\rm Let $T$ be a convex set. Denote by $\partial_e(T)$ the set of extremal points of
$T.$

\vspace{0.1in}

{\rm (18)} Let $X$ be a compact metric space with infinitely many points and $\af: X\to X$ be a
minimal homeomorphism. A Borel set $Y\subset X$ is said to be {\it universally null}, if
$
\mu(Y)=0
$
for all $\mu\in T_\af.$

\vspace{0.1in}

{\rm (19)} Let $A_\af$ be the simple crossed product. A crucial assumption that we make in this paper
is that
$\rho(K_0(A_\af))$ (see (4) above) is dense in
$Aff(T(A_\af)).$ }
\end{NN}

 We will use the following theorem (\cite{LP}).

\begin{thm}\label{LP}
Let $X$ be a finite dimensional compact metric space with infinitely many points and $\af: X\to X$ be
a minimal homeomorphism.  Then $A_\af$ has tracial rank zero if and only if
$\rho(K_0(A_\af))$ is dense in $Aff(T(A_\af)).$
\end{thm}

Minimal dynamical systems whose crossed product \CA s satisfy the above condition have be given and
discussed in \cite{LP}. It should be mentioned that if $(X,\af)$ is a minimal Cantor system, then
condition in \ref{LP} is always satisfied.


\section{Uniform Rohlin Tower Lemma  and  mean dimension zero}



\begin{Def}\label{U00}
{\rm Let $X$ be a compact metric space and let $\af: X\to X$ be a homeomorphism.
We say that $(X,\af)$ has the {\it small-boundary property} if for every point
$x\in X$ and every open neighborhood of $x$ there exists an open neighborhood
$V\subset U$ such that $\mu({\overline{V}}\setminus V)=0$ for all $\mu\in T_\af.$

By a result of Lindenstrauss and Weiss (see  \cite{LW}, \S 5), if $(X,\af)$  has small
boundary property then $(X,\af)$ has mean dimension zero. The converse
is also true, for example, if $(X,\af)$ is minimal (see Theorem 6.2 of \cite{Ld}).

It is also shown in \cite{LW} that if $X$ has finite covering dimension then any
minimal system $(X,\af)$ has mean dimension zero.

The following is an easy lemma.}
\end{Def}

\begin{Lem}\label{U0}
Let $X$ be a compact metric space with infinitely many points and let
$\af:X\to X$ be a homeomorphism. Suppose that $\partial_e(T_\af)$ is countable.
Then $(X,\af)$ has small boundary property. Consequently $(X,\af)$
has mean dimension zero.

More precisely, given any point $x\in X$ and $\dt>0,$ there is an open ball of $X$ with center at $x$
and radius $\dt/2<r<\dt$ such that
$$
\mu(\{y\in X: \di(x, y)=r\})=0
$$
for all $\mu\in T_\af.$

\end{Lem}
\begin{proof}
Let $\partial_e(T_\af)=\{\mu_1,\mu_2,...,\mu_n,...\}.$ Given a point
$x\in X$ and $\dt/2<r<\dt$
define
$$
R=\{y\in X: \dt/2<\di(y,x)<\dt\}\andeqn
$$
$$
 C_r=\{y\in X: \di(y,x)=r\}.
 $$
 Since
 $$
 \mu(R)=\mu(\cup_{\dt/2<r<\dt} C_r)
 $$
 and $\mu(R)\le 1$ for all $\mu\in T_\af,$ there are uncountably many $r\in (\dt/2, \dt)$
 such that
 $$
 \mu_n(C_r)=0,\,\,\,n=1,2,...,.
 $$
 Let $r$ be one of them. It follows that
$$
\mu(C_r)=0
$$
for all $\mu\in T_\af.$

\end{proof}

The  Rohlin Tower Lemma is well known in ergodic  theory. The following two lemmas are
 some uniform versions of it
which will be used later.

\begin{Lem}\label{Rok0}
Let $X$ be a compact metric space with infinitely many points, let
$\alpha: X\to X$ be a minimal homeomorphism and let $T_\af$ be
the set of $\alpha$ invariant probability measure.  Suppose that $(X,\af)$ has mean dimension zero.
Then, for any integer $n\ge 1,$ there exist finitely many open subsets
$G_1,G_2,..,G_m\subset X$
such that
$\alpha^j(G_i)$ are mutually disjoint for $0\le j\le h(i)-1,$ $0\le i\le m,$
$h(i)\ge n$ for each $i$ and
$\mu(X\setminus \cup_{i=1}^m \cup_{j=0}^{h(i)-1}\af^j(G_i))=0$ for all $\mu\in T_\af.$

\end{Lem}

\begin{proof}
We start with a non-empty  open subset $\Om\subset X$ such that
$\af^j(\overline{\Om})$ are pairwise disjoint for $0\le j\le n-1.$ This is possible since
$\af$ is minimal. By \ref{U0} and \ref{U00}, we may assume
that $\mu(\partial(\Om))=0$ for all $\mu\in T_\af.$

Let $Y=\overline{\Om}.$
For each $y\in Y,$ define
$$
r(y)=\min\{m>1: \af^m(y)\in Y\}.
$$

It follows from Theorem 2.3 of \ref{LP} (see also p.299 of \cite{LqP}) that
$\sup_{y\in Y} r(y)<\infty.$ Let $n(0)<n(1)<\cdots n(l)$ be distinct values in the range
of $r,$ and for each $0\le k\le l,$ set
$$
Y_k=\overline{\{y\in Y: r(y)=n(k)\}}\andeqn Y_k^o={\rm int}\{ y\in Y: r(y)=n(k)\}.
$$
Set
$$
X_k=\{y\in Y: r(y)\le n(k)\}.
$$
Since $Y$ is closed, so is  $X_k.$  Moreover,
$Y_0=X_0.$  Then
$$
Y_0=X_0, Y_1={\overline{X_1\setminus X_0}}, ..., Y_l=\overline{X_l\setminus X_{l-1}}.
$$

Note that $n(0)\ge n.$

 Denote $\Om_0={\rm int}(Y).$
Note that $\Om\subset \Om_0.$ Therefore $\overline{\Om_0}=Y.$ Put
$$
S_1=\af^{n(0)}(\Om_0)\cap \Om_0.
$$
Then $S_1$ is open and
\begin{align}
(\af^{n(0)}(Y)\cap Y)\setminus S_1&=
[(\af^{n(0)}(Y)\cap Y)\setminus  \af^{n(0)}(\Om_0)]\bigcup [(\af^{n(0)}(Y)\cap Y)\setminus \Om_0]\\
&\subset \af^{n(0)}(\partial(\Om_0))\bigcup\partial(\Om_0).
\end{align}
It follows that
$$
\mu((\af^{n(0)}(Y)\cap Y)\setminus S_1)=0
$$
 for all $\mu\in T_\af.$
 Note that $\af^{-n(0)}(\af^{n(0)}(Y)\cap Y))=Y_0.$
By the continuity of $\af,$ we also have
$$
 \af^{-n(0)}(S_1)=Y_0^o.
$$
It follows that
\beq\label{errok1}
\mu(X_0\setminus {\rm int}X_0)=\mu(Y_0\setminus Y_0^o ))=0
\eneq
for all $\mu\in T_\af.$ For $k>0,$ let
$$
S_k=\af^{n(k)}(\Om_0)\cap \Om_0.
$$
Then $S_k$ is open and, as above,
\beq\label{errok2}
\mu((\af^{n(k)}(Y)\cap Y)\setminus S_k))=0
\eneq
for all $\mu\in T_\af$ and $1\le k\le l.$  We have
\beq\label{errok2+}
 \af^{-n(k)}(\af^{n(k)}(Y)\cap Y)\setminus X_{k-1}=X_k\setminus X_{k-1}\andeqn \af^{-n(k)}(S_k)\setminus
 X_{k-1}=Y_k^o.
\eneq

Moreover, for $k>0,$ by (\ref{errok2+}),
\begin{align}
&X_k\setminus {\rm int}(X_k) \subset [(X_k\setminus X_{k-1})\setminus Y_k^o]
\bigcup(X_{k-1}\setminus {\rm int}(X_{k-1}))\\
&\subset (\af^{-n(k)}(\af^{n(k)}(Y)\cap Y)\setminus X_{k-1})\setminus (\af^{-n(k)}(S_k)\setminus
X_{k-1}))
\bigcup (X_{k-1}\setminus {\rm int}(X_{k-1})\\
&\subset  (\af^{-n(k)}(\af^{n(k)}(Y)\cap Y)\setminus \af^{-n(k)}(S_k))\bigcup (X_{k-1}\setminus {\rm
int}(X_{k-1})).
\end{align}
By induction on $k,$ combing the above with (\ref{errok1}) and with (\ref{errok2}), we conclude that
\beq\label{errok3}
\mu(X_k\setminus {\rm int}(X_k))=0
\eneq
for all $\mu\in T_\af,$ $1\le k\le l.$

We also have
\begin{align}
Y_k\setminus Y_k^o &\subset \overline{X_k\setminus X_{k-1}}\setminus (\af^{-n(k)}(S_k)\setminus X_{k-1})\\
&\subset (\af^{-n(k)}(\af^{n(k)}(Y)\cap Y)\setminus {\rm int}(X_{k-1}))\setminus
(\af^{-n(k)}(S_k)\setminus X_{k-1})\\
&\subset (\af^{-n(k)}(\af^{n(k)}(Y)\cap Y)\setminus \af^{-n(k)}(S_k))
\bigcup (X_{k-1}\setminus {\rm int}(X_{k-1})).
\end{align}
From this, by (\ref{errok2}) and (\ref{errok3}), we have
\beq\label{errok4}
\mu(Y_k\setminus Y_k^o)=0\,\rforal \mu\in T_\af.
\eneq

It follows from Theorem 2.3 of \cite{LP} (see also p.299 of \cite{LqP}) that

(i) $\af^j(Y_k^o)$ are pairwise disjoint for $1\le j\le n(k),$ $0\le k\le l;$

(ii) $\bigcup_{k=0}^l\bigcup_{j=0}^n(k)\af^j(Y_k)=X.$

Moreover
$$
\mu(X\setminus \bigcup_{k=0}^l\bigcup_{j=0}^{n(k)}\af^j(Y_k^o))
\le \sum_{k=0}^l\sum_{j=0}^{n(k)}\mu(\af^j(Y_k\setminus Y_k^o))=0
$$
for all $\mu\in T_\af.$ Define
$G_k=\af(Y_k^o),$ $k=0,1,...,l.$ With $m=l+1$ and $h(k)=n(k)+1,$ we see the lemma follows.

\end{proof}

\begin{Lem}\label{Rok1}
Let $X$ be a compact metric space with infinitely many points, let
$\alpha: X\to X$ be a minimal homeomorphism and let $T_\af$ be
the set of $\alpha$ invariant probability measure.  Suppose that $(X,\af)$ has mean dimension zero.
Let $\{y_1,y_2,...,y_k\}$ be $\eta/3$-dense subset of $X$ for some $\eta>0.$

Then, for any integer $n\ge 1,$ there exists an open subset $G\subset X$ containing a subset
$\{x_1,x_2,...,x_k\}$ which is $\eta$-dense in $X$
with $\di(x_i,y_i)<\eta/3$ ($1\le i\le k$)
such that
$\alpha^i(G)$ are mutually disjoint for $0\le i\le n-1$ and
$\mu(\cup_{i=0}^{n-1}\alpha^i(G))>1-\ep$ for all $\mu\in T_\af.$

Moreover,
$$
\mu(\partial(G))=0
$$
for all $\mu\in T_\af.$


\end{Lem}

\begin{proof}
Choose an integer $K>0$ such that $1/K<\ep.$ Let $N=nK.$ By \ref{Rok0}, we obtain finitely many open
subsets $G_1,G_2,...,G_m$ such that

(i) $\af^j(G_i)$ are pairwise disjoint for $1\le i\le m,$ $0\le j\le h(i);$

(ii) $h(i)\ge N,$ $1\le i\le m;$

(iii) $\mu(X\setminus \cup_{i=1}^m\cup_{j=0}^{h(i)-1}\af^j(G_i))=0$ for all $\mu\in T_\af.$

Write $h(i)=L(i)n+r(i),$ where $L(i)\ge 1$ and $n>r(i)\ge 0$ are integers, $i=1,2,...,m.$ Define, for
each $i,$
$$
U(i,1)=\af^n(G_i), U(i,2)=\af^{2n}(G_i),...,U(i,L(i)-1)=\af^{(L(i)-1)n}(G_i).
$$
Note that
\beq\label{rrok1m}
\mu(G_i)\le  {1\over{nK}}\mu(\cup_{j=0}^{h(i)-1}\af^j(G_i)),\,\,\,1\le i\le m
\eneq
for all $\mu\in T_\af.$

So
\beq\label{rrok1m2}
\mu(\cup_{j=L(i)}^{h(i)-1}\af^j(G_i))=r(i)\mu(G_i)\le {1\over{K}}\mu(\cup_{j=0}^{h(i)-1}\af^j(G_i))
\eneq
for all $\mu\in T_\af$ and $1\le i\le m.$

Let $G=\cup_{i=1}^m G_i\bigcup (\cup_{i=1}^m \cup_{s=1}^{L(i)-1} U(i,s)).$ Then

(1) $\af^j(G)$ are pairwise disjoint for $0\le j\le n-1,$

and, by (iii) and by (\ref{rrok1m2}),

(2) $\mu(\cup_{j=0}^{n-1}\af^j(G))>1-\sum_{i=1}^m \mu(\cup_{j=L(i)}^{h(i)-1}\af^j(G_i))
>1-{1\over{K}}>1-\ep$ for
all $\mu\in T_\af.$

Now
let $\{y_1,y_2,...,y_k\}$ be an $\eta/3$-dense set.
Define $y_i'=\af^{-1}(y_i),$ $i=1,2,...,k.$
Choose $\dt>0$ such that
$$
\di(\af(x), \af(y))<\eta/9
$$
whenever $\di(x,y)<\dt.$

Choose
$z_1=y_1'.$ Since $y_2'$ is a condensed point, there is $z_2\in
O(y_2'),$ where $O(y_2')=\{x\in X: \di(y_2,x)<\dt\},$ such that
$z_2\not\in \{\alpha^n(x_1): n\in {\mathbb Z}\}.$ We then choose
$z_3\not\in \{\alpha^n(x_1), \alpha^n(x_2): n\in {\mathbb Z}\}$
but $\di(z_3, y_2)<\dt.$ By induction, we obtain
$\{z_1,z_2,...,z_k\}\subset X$ such that none of $z_i$ lies in the
obit of $z_j$ if $i\not=j.$ We note that $\{\af(z_1),\af(z_2),...,\af(z_k)\}$ is
$4\eta/9$-dense in $X.$
So we may start with an open subset $\Om$ which contains $\{z_1,z_2,...,z_k\}$
at the beginning of the proof of \ref{Rok0}.

Note that, by the proof of \ref{Rok0},
 $G_k=\af(Y_k^o),$ $k=0,1,...,l.$  In the proof of \ref{Rok0},
$$
\bigcup_{k=0}^lY_k\supset Y={\overline{\Om}}.
$$
It follows that
$$
\af(Y)\setminus \bigcup_{k=0}^l G_k\subset \bigcup_{k=0}^l\af(Y_k\setminus Y_k^o).
$$
Since
$$
\mu(Y_k\setminus Y_k^o)=0
$$
for all $\mu\in T_\af,$ and since $\af$ is minimal, for each $i,$
$$
U(\af(z_i))\cap \cup_{k=1}^mG_k \not=\emptyset,
$$
where $U(\af(z_i))=\{x\in X: \di(\af(z_i),x)<\eta/9\}.$
Choose a point $x_i\in U(\af(z_i))\cap \cup_{k=1}^lG_k,$ $1\le i\le k.$
Then the above proof shows that
$$
x_i\in G,\,\,\, i=1,2,...,k.
$$
Note that $\di(x_i, y_i)<\eta/3$ $i=1,2,...,k$ and $\{x_1,x_2,...,x_k\}$ is
$\eta$-dense in $X.$

\end{proof}

\vspace{0.2in}

Let $X$ be a compact metric space and let $A$ be a unital \CA. Suppose that $\phi: C(X)\to A$ is a
\hm. Then $\phi$ can be extended to a \hm\, from ${\cal B}(X),$ the algebra of all bounded Borel
functions, to the enveloping von-Neumann algebra $A^{**}.$

\begin{Lem}\label{Lsemi2}
Let $X$ be a compact metric space and $\phi: C(X)\to A$ be a unital monomorphism from $C(X)$ into a
unital simple \CA\, $A.$ Suppose that
$$
\mu_{\tau}(\overline{G}\setminus G)=0
$$
for all $\tau\in T(A).$

Then $\phi(\chi_{G})$  (in $A^{**}$) is continuous function on
$T(A),$ or equivalently, for any $\ep>0,$ there exists $f\in
C(X),$ with $0\le f(t)\le 1$ for all $f\in X$ and $f(t)=0$ if
$t\in X\setminus G$ such that
$$
|\tau(\phi(f))-\mu_{\tau}(G)|<\ep
$$
for all $\tau\in T(A).$
\end{Lem}

\begin{proof}
Note that $X\setminus {\bar G}$ is open. There are increasing sequences of nonnegative functions
$\{g_n\}, \{h_n\}\subset  C(X)$ with $0\le g_n\le \chi_G$ and $0\le h_n\le 1-\chi_{{\bar G}}$ such
that
$$
\lim_{n\to\infty}g_n(t)=\chi_G(t) \andeqn  \lim_{n\to\infty} h_n(t)=1-\chi_{\bar G}
$$
for each $t\in X.$ Let
$$
I_1=\{f\in C(X): f(t)=0\rforal t\in X\setminus G\} \andeqn I_2=\{f\in C(X): f(t)=0\rforal f\in
\bar{G}\}.
$$
Let $B_1={\overline{\phi(I_1)A\phi(I_1)}}$ and
$B_2={\overline{\phi(I_2)A\phi(I_2)}}$ be the hereditary \SCA s
corresponding to the open subsets $G$ and $X\setminus {\bar G},$ respectively. It follows easily
that
$$
\lim_{n\to\infty}\tau(g_n)=\tau(\chi_G)\andeqn \lim_{n\to\infty}\tau(h_n)=\tau(1-\chi_{\bar G})
$$
for all $\tau\in T(A).$ Note that $\tau((1-\chi_{\bar G})+\chi_G))=1$ for all $\tau\in T(A).$ Consider the
function
$\widehat{\phi(g_n+h_n)}$ on $T(A)$ defined by
$\widehat{\phi(g_n+h_n)}(\tau)=\tau(\phi(g_n+h_n)$ for $\tau\in
T(A).$ By the Dini  Theorem, $\phi(g_n+h_n)$ converges uniformly to 1 on $T(A).$ In other words, for
any $\ep>0,$ there exists
$N>0$ such that
$$
|\tau(\chi_G-g_n)+\tau((1-\chi_{\bar G})-h_n)| <\ep
$$
for all $\tau\in T(A)$ provided that $n\ge N.$ Since both
$\tau(\phi(\chi_G-g_n))$ and $\tau(\phi((1-\chi_{\bar G})-h_n))$
are positive, one has, if $n\ge N,$
$$
\tau(\phi(\chi_G-g_n))<\ep
$$
for all $\tau\in T(A).$ This shows that $\tau(\phi(g_n))$ converges to $\tau(\phi(\chi_G))$ uniformly
on $T(A).$ It follows that $\phi(\chi_G)$ is a continuous function on $T(A).$
\end{proof}

\section{Perturbations}

{

The following lemma is well-known (note that finite dimensional \CA s are semiprojective and their
unit balls are compact).

\begin{Lem}\label{pertM}
Let $F$ be a finite dimensional \CA. Then for any $\ep>0$
there exist a finite subset ${\cal G}\subset F$ and $\dt>0$ satisfying the following: For any ${\cal
G}$-$\dt$-multiplicative \morp\, $\phi: F\to A,$ where $A$ is any \CA, there exists a \hm\, $h: F\to A$
such that
$$
\|h-\phi|_F\|<\ep.
$$
\end{Lem}

\begin{Lem}{\rm (Lemma 4.1 of \cite{Lncd})}\label{Ltr}
Let $A$ be a unital \CA. For any $\ep>0$ and finite subset ${\cal F}\subset A,$ there exist a finite
subset ${\cal G}\subset A$ and
$\dt>0$ satisfying the following:

If $B$ is another unital \CA, $\phi: A\to B$ is a unital \morp\, which is ${\cal G}$-$\dt$-
multiplicative and $\tau\in T(B),$ then there exists a tracial state $t\in T(A)$ such that
$$
|\tau\circ \phi(a)-t(a)|<\ep
$$
for all $a\in {\cal F}.$
\end{Lem}

\begin{Lem}\label{Pert}
Let $X$ be a compact metric space with infinitely many points and
let $\alpha: X\to X$ be a minimal homeomorphism. Let
$G_1,G_2,...,G_L$ be finitely many open subsets with the property that
$\mu({\overline{G_i}}\setminus G_i)=0$ for all $\mu\in T_\af.$

For any $\ep>0$ and a finite subset ${\cal F}\subset C(X),$
 there exist a finite subset ${\cal G}_1\subset
C(X)$ and $\eta>0$ satisfying the following:

if there exists a projection $p\in A_{\alpha}$ and a unital \hm\,
$\phi_0: C(X)\to pAp$ with finite dimensional range such that

{\rm (1)} $\|pj_{\alpha}(f)-j_{\alpha}(f)p\|<\eta$ for all $f\in {\cal G}_1,$

{\rm (2)} $\|pj_{\alpha}(f)p-\phi_0(f)\|<\eta$ for all $f\in {\cal G}_1$ and

{\rm (3)} $\tau(1-p)<\eta$ for all $ \tau\in T(A_\af),$

\noindent and if $\phi: A_{\alpha}\to M_k$ is a unital ${\cal G}_2$-$\dt$-multiplicative \morp\, (for
some $k>0$) ,
where ${\cal G}_2$ is a finite subset of $A_\af$ and $\dt>0,$ both of which depend  on the image of
$\phi_0,$ ${\cal G}_1,$ $\eta,$ $\ep$ as well as  $G_1,G_2,...,G_L,$

\noindent then there is $\tau\in T(A_{\alpha}),$ such that
$$
|{\rm tr}\circ \phi\circ j(g)-\tau(g)|<\ep/2\andeqn |{\rm tr}\circ \phi\circ\phi_0(g)-\tau(g)|<\ep
$$
for all $g\in {\cal G}_1,$ there are $\{y_1,y_2,...,y_K\}\subset X$ and mutually orthogonal rank one
projections in $M_k$ such that
$$
\|\sum_{i=1}^Kf(y_i)p_i-\phi\circ\phi_0\circ (f)\|<\ep
$$
for all $f\in {\cal F}$
and
$$
\mu_{\tau}(G_j)+\ep>{N_j\over{k}}> \mu_{\tau}(G_j)-\ep,
$$
where $N_j$ is the number of $y_i's$ in $G_j.$
Moreover, ${k-K\over{k}}<\ep.$
\end{Lem}

\begin{proof}
To simplify notation, without loss of generality, we may assume that ${\cal F}$ is in the unit ball of
$C(X).$




 Let
$$
\gamma_0=\inf\{\mu_{\tau}(G_j): \mu\in T(A), j=1,2,...,L\}.
$$
Since $A_{\alpha}$ is simple, one has $\gamma_0>0.$
 By Lemma \ref{Lsemi2}, choose
$g_j\in C(X)$ with $0\le g_j\le 1,$ $g_j(x)=0$ if $x\not\in G_j$
and
\begin{eqnarray}\label{eppert1}
\mu_{\tau}(G_j)<\tau(j_{\alpha}(g_j))-\min (\gm_0/4,\ep/8)
\end{eqnarray}
for all $\tau\in T(A)$ and $j=1,2,...,L.$

Let ${\cal F}_1={\cal F}
\cup\{g_j: 1\le j\le L\}.$ Let $\eta_1>0$ be such that
$$
|f(x)-f(x')|<\ep/4
$$
if $\di(x,x')<\eta_1$ for all $f\in {\cal F}_1.$
 Let $\eta=\min\{\gamma_0/32, \ep/64,\eta_1/32\}.$
Let ${\cal G}_1={\cal F}_1.$ Suppose that $p\in A_{\alpha}$ and
$\phi_0: C(X)\to pA_{\alpha}p$ is a \hm\, with finite dimensional
range which satisfies (1)-(3) as described in the statement (for the above ${\cal G}_1$ and $\eta$).

Put ${\cal F}_2=j_{\alpha}({\cal F}_1)\cup \phi_0({\cal F}_1)\cup\{p, 1-p\}\cup \{pj_{\alpha}(f)p:
f\in {\cal F}_1\}.$

 Let
${\cal G}\subset A_\af$ be a finite subset and $\dt>0$ be a positive number
required by Lemma \ref{Ltr} corresponding to ${\cal F}_2$  and
$\eta.$ Let $C$ be the image of $\phi_0$  which is a finite dimensional \CA.
Choose a smaller $\dt$ required by \ref{pertM} and a larger ${\cal G}$ which contains a finite subset
required by \ref{pertM} for $C$ and $\eta.$


Let ${\cal G}_2={\cal G}\cup {\cal F}_2.$
Now let $\phi: A_{\alpha}$ be a unital ${\cal G}_2$-$\dt$-multiplicative \morp\, from $A_{\alpha}\to
M_k$ (for some $k>0$).

By \ref{pertM} (and choice of ${\cal G}$ and $\dt$),
 we may also assume that there is a \hm, $\phi_{00}:
C(X)\to EM_kE$ (for some projection $E$) such that
$$
\|\phi_{00}(f)-\phi\circ \phi_0(f)\|<\eta
$$
for all $f\in {\cal F}_1.$

By the choice of ${\cal G}$ and $\dt,$ applying \ref{Ltr}, there is a tracial state
$\tau\in T(A)$ such that
$$
|\tau(a)-{\rm tr}\circ \phi(a)|<\eta
$$
for all $f\in {\cal F}_2.$ In particular,
$$
|\tau(1-p)-{\rm tr}\circ \phi(1-p)|<\eta.
$$
It follows that
\begin{eqnarray}\label{eTr(1-p)}
{\rm tr}\circ\phi(1-p)<2\eta<\ep/4.
\end{eqnarray}

Moreover
$$
|\tau(j_{\alpha}(f))-{\rm tr}\circ \phi_{00}(f)|<3\eta
$$
for all $f\in {\cal F}_1.$

Write $\phi_{00}(f)=\sum_{i=1}^Kf(y_i)p_i$ for all $f\in C(X),$ where $y_i\in X$ and
$\{p_1,p_2,...,p_K\}$ is a set of mutually orthogonal rank one projections in $M_k,$ and $0<K<k.$




On the other hand,
\begin{eqnarray}\label{eppert2}
|{\rm tr}(\phi_{00}(g_i))-\tau(j_{\alpha}(g_i))|<3\eta
\end{eqnarray}
for $i=1,2,...,L.$ It follows  from (\ref{eppert1}) and (\ref{eppert2}) that
$$
\mu_{\tau}(G_j)+\ep/2>{N_j\over{k}}>\mu_{\tau}(G_j)-\ep/2,
$$
 where $N_j$ is the number of $y_j$'s which lie
in $G_j,$ $j=1,2,...,L.$

By (\ref{eTr(1-p)}), we compute that
$$
{k-K\over{k}}<\ep/4<\ep.
$$


\end{proof}

\begin{Lem}\label{dig}
Let $X$ be a finite dimensional compact metric space with infinitely many points and
$\alpha: X\to X$ be a minimal homeomorphism.
Suppose that $\rho(K_0(A_\af))$ is dense in $Aff(T(A_\af).$ Then, any $\ep>0,$ $\sigma>0$ and finite
subset ${\cal F}\subset C(X),$ there are mutually orthogonal projections
$\{p_1,p_2,...,p_m\}\subset A_{\alpha}$ and $\{x_1,x_2,...,x_m\}\subset X$
such that

{\rm (1)} $\|pj_{\alpha}(f)-j_{\alpha}(f)p\|<\ep$ for $f\in {\cal F},$ where $p=\sum_{k=1}^mp_k,$

{\rm (2)} $\|pj_{\alpha}(f)p-\sum_{k=1}^m f(x_i)p_k\|<\ep$ for all $f\in {\cal F}$ and

{\rm (3)} $\tau(1-p)<\sigma$ for all $\tau\in T(A_{\alpha}).$


\end{Lem}

\begin{proof}
Let $\eta>0$ be such that
$$
|f(x)-f(x')|<\ep/4
$$
if $\di(x,x')<\eta$ for all $f\in {\cal F}.$ It follows from \ref{U00} that there are pairwise disjoint
open subsets $O_1,O_2,...,O_m$ with diameters less than $\eta$ such that
$\{\overline{O_1}, \overline{O_2},..., \overline{O_m}\}$ covers
$X$ and
\begin{eqnarray}\label{edig1}
\mu(\cup_{i=1}^m  (\overline{O_i}\setminus O_i))=0
\end{eqnarray}
for all $\mu\in T_\af.$

Let $g_i\in C(X)$ such that $0\le g_i\le 1,$ $g_i(x)=0$ if $x\not\in O_i$ and
$g(x)>0$ if $x\in O_i,$ $i=1,2,...,m.$
Let $B_i=\overline{j_\af(g_i)A_\af j_\af(g_i)}$ be the hereditary \SCA\, associated with the open set
$O_i,$ $i=1,2,...,m.$ Since $\rho(K_0(A_\af))$ is dense in $Aff(T(A_\af)),$ by \ref{LP},
$TR(A_\af)=0.$ In particular, $A_\af$ has real rank zero.
So $B_i$ contains an approximate identity consisting of projections
$\{e(i,n): n=1,2,...\},$ $i=1,2,...,m.$
It is easy to see that
$$
\lim_{n\to\infty}\tau(e(i,n))=\mu_{\tau}(O_i),
$$
where $\mu_{\tau}$ is the probability measure induced by the trace $\tau,$ for all $\tau\in T(A_\af).$
Put $e_n=\sum_{i=1}^me(i,n).$ By (\ref{edig1}),
$$
\lim_{n\to\infty}\tau(e_n)=1
$$
for all $\tau\in T(A_\af).$ Since $\{e_n\}$ is an increasing sequence, by the Dini Theorem,
$\{\widehat{e_n}\}$ converges uniformly on $T(A_\af).$ Therefore, for any $\sigma>0,$ there exists an
integer $n>0$ such that
\begin{eqnarray}\label{edig2}
\tau(e_n)>1-\sigma
\end{eqnarray}
for all $\tau\in T(A_\af).$ Define $p=1-e_n$ and $p_i=e(i,n),$ $i=1,2,...,m.$ Choose $x_i\in O_i.$ By
the choice of $\eta,$ one checks easily that (1), (2) and (3) follows (see the proof of 4.1 of
\cite{EGLP}).
\end{proof}

\begin{Lem}\label{Pert2}
Let $X$ be a finite dimensional compact metric space with infinitely many points and  let $\alpha:
X\to X$ be a minimal homeomorphism. Suppose that $\rho(K_0(A_\af))$ is dense in
$Aff(T(A_\af).$

Let
$G_1,G_2,...,G_L$ be finitely many open subsets with the property that
$\mu({\overline{G_i}}\setminus G_i)=0$ for all $\mu\in T_\af.$ For any $\ep>0$ and any finite subset
${\cal F}\subset C(X),$
there are a (specially selected)  projection $p\in A_{\af}$ with
$\tau(1-p)<\ep/2$ for all $\tau\in T(A_\af),$ and a finite subset
${\cal G}\subset A_{\alpha} $ and $\dt>0$ satisfying the
following:

if $\phi: A_{\alpha}\to M_k$ is a unital ${\cal G}$-$\dt$-multiplicative \morp\, (for some $k>0$),
then there is $\tau\in T(A_{\alpha})$  such that
$$
|{\rm tr}\circ \phi\circ j(g)-\tau(g)|<\ep/2\andeqn |{\rm tr}\circ \phi(pgp)-\tau(g)|<\ep
$$
for all $g\in {\cal F},$ and there are $\{y_1,y_2,...,y_K\}\subset X$ and mutually orthogonal rank one
projections $\{p_1,p_2,...,p_k\}$ in $M_k$ such that
$$
\|\sum_{i=1}^Kf(y_i)p_i-\phi\circ (pfp)\|<\ep
$$
for all $f\in {\cal F}$
 and
$$
\mu(G_j)+\ep>{N_j\over{k}}> \mu(G_j)-\ep,
$$
where $N_j$ is the number of $y_i's$ in $G_j$ and $\mu$ is the probability measure induced by $\tau.$
Moreover, ${k-K\over{k}}<\ep.$
\end{Lem}

\begin{proof}
To prove this lemma, we combine \ref{Pert} and \ref{dig}. Fix $\ep>0$ and a finite subset ${\cal
F}\subset C(X).$ Let ${\cal G}_1\subset C(X)$ be a finite subset and $\eta>0$ required by \ref{Pert}.
By applying \ref{dig}, we obtain a projection $p\in A_\af$ and a unital \hm\, $\phi_0: C(X)\to pAp$
with finite dimensional range which satisfies (1)-(3) in \ref{Pert}. We then apply \ref{Pert} to
obtain this lemma.

\end{proof}
\begin{Lem}\label{LMM}
Let $A$ be a unital  simple \CA\, with the property: two
projections $p$ and $q$ in $A$ with $\tau(p)=\tau(q)$ for all
$\tau\in (A)$ are equivalent.

Let $X$ be a compact metric space and $h_1, h_2: C(X)\to A$ be two
unital monomorphisms. Suppose that
\begin{eqnarray}\label{eLMM}
\tau\circ h_1(f)=\tau\circ h_2(f)
\end{eqnarray}
for all $\tau\in T(A).$ Suppose also that, for any $r>0,$ there
are finitely many pairwise disjoint open subsets $U_1,U_2,...,U_m$
whose diameters are less than $r$  such that $X=\cup_{i=1}^m
{\overline{U_i}}$ and
$$
\mu_{\tau\circ h_1}(\cup_{i=1}^m(\overline{U_i}\setminus U_i))=0
$$
for all $\tau\in T(A).$

 Then,  for any $\eta>0,$  there exist a
finite subset ${\cal F}_0\subset C(X),$ ${\cal F}\subset A$ and
$\dt>0$ satisfying the following: for any ${\cal
F}$-$\dt$-multiplicative \morp\, $\phi: A\to B$ and any \hm\,
$\psi_1,\psi_2: C(X)\to B$ for some unital stably finite \CA\, $B$ with
$$
\|\phi\circ h_i(f)-\phi_i(f)\|<\dt
$$
for all $f\in {\cal F}_0,$ $i=1,2,$ one has
$$
\mu_{t\circ \psi_1}(S)\le  \mu_{t\circ  \psi_2}(B_\eta(S))\andeqn
\mu_{t\circ \psi_2}(S)\le \mu_{t\circ \psi_1}(B_\eta(S))
$$
for any $t\in T(B)$ and any closed subset $S\subset X,$
where
$B_\eta(S)=\{x\in X: \di(x,S)<\eta\}.$
\end{Lem}

\begin{proof}
Fix $\eta>0.$ Let $X=\sum_{j=1}^NX_i,$ where each $X_i$ is a
clopen set which is $\eta/4$-connected, i.e., for any
two points $x, y\in X_i,$ there are $x_1,x_2,...,x_m\in X_i$
such that $\di(x,x_1)<\eta/4,$ $\di(x_i,x_{i+1})<\eta/4$ and
$\di(x_m,y)<\eta/4.$

Let $U_1,U_2,...,U_m$ be pairwise disjoint non-empty open subsets
whose diameters are less than $\eta/8$ such that $X=\cup_{i=1}^m
\overline{U_i}$ and
$$
\mu_{\tau\circ h_1}(\cup_{i=1}^m(\overline{U_i}\setminus U_i))=0
$$
for all $\tau\in T(A).$

Let
$$
d=\inf\{\mu_{\tau\circ h_1}(U_i): 1\le i\le m, \tau\in T(A)\}.
$$
Since $A$ is simple, $d>0.$

Let $e_1=h_1(\chi_{X_i})$ and $f_i=h_2(\chi_{X_i}),$ where
$\chi_{X_i}$ is the characteristic function on the clopen set
$X_i,$ $i=1,2,...,N.$ Then, for any $\tau\in T(A),$
\begin{eqnarray}\label{eLMM1}
\tau(e_i)=\tau(f_i)
\end{eqnarray}
for all $\tau\in T(A).$ By the assumption on $A,$ there is a
partial isometry $u_i\in A$ such that
\begin{eqnarray}\label{eLMM2}
u_i^*u_i=e_i\andeqn u_iu_i^*=f_i\,\,\,i=1,2,...,N.
\end{eqnarray}

Let $\Lambda$ be a subset of $\{1,2,...,m\}.$
By \ref{Lsemi2}, for each $\Lambda,$ there exists a $g_{\Lambda}\in C(X)$ with
$0\le g_{\Lambda}\le 1,$ $g_{\Lambda}(x)=1$ if $x\in \cup_{i\in \Lambda}U_i$ and $g_i(x)=0$ if
$\di(x,\cup_{i\in\Lambda}U_i)>\eta/128$
such that
\begin{eqnarray}\label{eLMM3}
\tau(h_1(g_{\Lambda}))-{d\over{8}}<\mu_{\tau\circ h_1}(\cup_{i\in \Lambda}U_i)
\end{eqnarray}
 for all $\tau\in T(A),$ $i=1,2,...,m.$

Let ${\cal F}_0=\{g_{\Lambda}: \Lambda\subset \{1,2,...,m\}\}$ and ${\cal F}=\{u_i,u_i^*:
1\le i\le N\}\cup_{i=1}^2h_i({\cal F}_0).$
 Let ${\cal G}$ be a finite subset and $\dt>0$ be required
by \ref{Ltr} for the above $A,$ ${\cal F}$ and $d/8.$ We may assume
that $\dt<d/4.$

Now suppose that $\phi: A\to B$ is a ${\cal
G}$-$\dt/4$-multiplicative \morp\, and $\psi_i: C(X)\to B$ is (for
each $i$) a \hm\, such that
\begin{eqnarray}\label{eLMM4}
\|\psi_i(f)-\phi\circ h_i(f)\|<\dt/4
\end{eqnarray}
for all $f\in {\cal F}.$

Hence
\begin{eqnarray}\label{eLMM5}
\|\psi_1(\chi_{X_i})-\phi_1(u_i)\phi(u_i)^*\|<\dt \andeqn
\|\psi_2(\chi_{X_i})-\phi(u_i)^*\phi(u_i)\|<\dt
\end{eqnarray}
$i=1,2,...,N.$ With $\dt<d/4<1,$ it follows (for example, from 2.5.3 of \cite{Lnb})
that $\psi_1(\chi_{X_i})$ is equivalent to $\psi_2(\chi_{X_i})$ in $B,$
$i=1,2,...,N.$

In particular,
\begin{eqnarray}\label{eLMM6}
t(\psi_1(\chi_{X_i}))=t(\psi_2(\chi_{X_i}))
\end{eqnarray}
for all $t\in T(B),$ $i=1,2,...,N.$

By the choice of ${\cal G}$ and $\dt,$ applying \ref{Ltr},
we have,
for each $t\in T(B),$ there is $\tau\in T(A)$ such
that
\begin{eqnarray}\label{eLMM7}
|\tau(h_1(g_{\Lambda})-t\circ\psi_j(g_{\Lambda})|<d/8
\end{eqnarray}
$j=1,2$ and $\Lambda\subset \{1,2,...,m\}.$

For any closed subset $S\subset X,$ if $S$ is a union of some of
$X_i,$ then, by (\ref{eLMM6}),
\begin{eqnarray}\label{eLMM8}
\mu_{t\circ\psi_1}(S)=\mu_{t\circ \psi_2}(S).
\end{eqnarray}

Suppose that $S$ is a closed subset of $X$ which is not a finite
union of some $X_i$'s. Then, there must be a point $\xi\in
B_{5\eta/16}(S)\setminus B_{\eta/4}(S).$ But $\di(\xi, U_j)=0$ for
some $j.$ Since diameter of $U_j$ is less than $\eta/8,$
\begin{eqnarray}\label{eLMM9}
U_j\subset B_{7\eta/16}(S)\subset B_{\eta/2}(S).
\end{eqnarray}
It follows from (\ref{eLMM7} that
\begin{eqnarray}\label{eLMM10}
\mu_{t\circ \psi_i}(U_j)>d/2
\end{eqnarray}
for all $t\in T(B),$ $i=1,2.$
Since $U_j\cap B_{7\eta/64}(S)=\emptyset,$
we have,
\begin{eqnarray}\label{eLMM11}
\mu_{t\circ \psi_i}(B_{\eta}(S))>d/2+\mu_{t\circ
\psi_i}(B_{7\eta/64}(S))
\end{eqnarray}

There is a $\Lambda\subset \{1,2,...,N\}$ such that
$\cup_{i\in \Lambda}{\overline{U_i}}\supset S.$ Suppose that $\Lambda$ is smallest
such subset of $\{1,2,...,N\}.$
Then
\begin{eqnarray}\label{eLMM12}
{\rm supp}g_{\Lambda} \subset B_{7\eta/64}(S)
\andeqn
\mu_{t\circ\psi_i}(B_{7\eta/64}(S))\ge t(\psi_i(g_{\Lambda}))
\end{eqnarray}
for all $t\in T(B)$ and $i=1,2.$

By \ref{eLMM7},
\begin{eqnarray}\label{eLMM13}
|t\circ \psi_1(g_{\Lambda})-t\circ \psi_2(g_{\Lambda})|<d/8
\end{eqnarray}
for all $t\in T(B).$
It follows that, by applying (\ref{eLMM13}), (\ref{eLMM12}) and (\ref{eLMM11}),
\begin{eqnarray}\label{eLMM14}
\mu_{t\circ \psi_1}(S)\le t(\psi_1(g_{\Lambda}))\le t(\psi_2(g_{\Lambda}))+d/8
\le \mu_{t\circ\psi_2}(B_{7\eta/64}(S))+d/8\le \mu_{t\circ\psi_2}(B_{\eta})
\end{eqnarray}
for all $t\in T(B).$
Similarly,
\begin{eqnarray}
\mu_{t\circ \psi_2}(S)\le  \mu_{t\circ\psi_1}(B_{\eta})
\end{eqnarray}
for all $t\in T(B).$

\end{proof}

\begin{Lem}\label{Permutation}
Let $X$ be a  finite dimensional compact metric space with infinitely many points and let $\alpha:
X\to X$ be a minimal homeomorphism. Suppose that
$\rho(K_0(A_\af))$ is dense in
$Aff(T(A_\af)).$

 Let  $\ep>0$ and let ${\cal F}\subset C(X)$ be a finite
subset. Let $\eta>0$ be any positive number such that
$$
|f(t)-f(t')|<\ep/8
$$
if $\di(t,t')<\eta$ for all $f\in {\cal F}.$

Let $n$ be an integer so that $1/n<\ep/4$ and let
$G$ be an open set
such that $\af^j(G)$ are  pairwise
disjoint for $0\le j\le n-1$  with the following properties:

{\rm (i)} $G$ contains  $x_i,$ $i=1,2,...,l,$ where $\{x_1,x_2,...,x_l\}$ is $\eta/2$-dense in $X,$




{\rm (ii)}
$\mu(\cup_{j}\af^j(G))>1-\ep/16$ for all $\mu\in T_\af$ and

{\rm (iii)} $\mu(\partial(G))=0$ for all $\mu\in T_\af.$

Then there exist a (specially selected) projection $p\in A_\af$ with
$\tau(1-p)<\ep/2$ for all $\tau\in T(A_\af),$
a finite subset  ${\cal G}\subset A_{\alpha}$ and
$\dt>0$
satisfying the following:

if $\phi: A_{\alpha}\to M_k$ (with $k>ln$) is a ${\cal G}$-$\dt$-multiplicative \morp,
then there are  $m$ distinct points
$$
\{y_i, i=1,2,...,m\}
$$
with $y_i\in G,$ $x_i=y_i,$ $i=1,2,...,l\le m$   and ${k-mn\over{k}}<\ep/4$ such  that
\begin{eqnarray}\label{ePms}
\|\sum_{j=0}^{n-1}\sum_{i=1}^mf(\alpha^j(y_i))p_{i,j}+\sum_{i=K+1}^Nf(z_i)p_i-\phi
(pj_{\alpha}(f)p)\|<\ep
\end{eqnarray}
($K=mn<N<k$) for all $f\in {\cal F},$ where
$$
\{p_{i,j}: 1\le i\le m, 0\le j\le n-1 \}\cup\{p_{K+1},...,p_N\}
$$ is a set
of mutually orthogonal rank one projections in $M_k$  and
$\{z_{K+1},...,z_N\}\subset X.$
\end{Lem}

\begin{proof}
Let $\eta_1>0$ such that $\eta_1<\eta$ and
\begin{eqnarray}\label{epPert1}
\di(\af^j(x),\af^j(x'))<\eta/2
\end{eqnarray}
if $\di(x,x')<\eta_1,$ $-n+1\le j\le n-1.$
Let $\eta_2>0$ be such that $\eta_2<\eta_1$ and
\begin{eqnarray}\label{epPert1+}
\di(\af^j(x),\af^j(x'))<\eta_1/2
\end{eqnarray}
if $\di(x,x')<\eta_2,$ $j=1,2,...,n-1.$









Since $X$ has finite covering dimension, $(X,\af)$ has mean dimension zero (see \ref{U00}).
Let $U_i$ be an open ball with center at $x_i$ and radius $\eta_2/4$ such that
$\mu({\overline{U_i}}\setminus U_i)=0$ for all $\mu\in T_\af,$ $i=1,2,...,L.$

Now we apply \ref{Pert2} with open subsets $\{U_i: 1\le i\le L\}$ and
$\{\af^j(G): 0\le j\le n-1\}.$
Let $\dt_1>0.$
By \ref{Pert2} for ${\ep\over{8(n+1)}}$  and ${\cal F},$ with sufficiently large
${\cal G}$
and sufficiently small
$\dt,$ we may assume that $k$ is sufficiently large and
\begin{eqnarray}\label{eppm1}
\|\phi\circ (pj_{\alpha}(f)p)-\sum_{i=1}^Nf(z_i)p_i\|<\min
\{\ep/8,\dt_1 \}
\end{eqnarray}
where $p\in A_{\alpha}$ is a specially selected projection with
$\tau(1-p)<\ep/8$ for all $\tau\in T(A_{\alpha}),$ where
${k-N\over{k}}<\ep/8$ and where $\{z_1,...,z_N\}$ is a set of
distinct points of $X.$
By applying \ref{Pert2} (with finitely many open sets $U_i$'s and
$\af^j(G)$'s in place of $G_i$) and using (ii) above,
we may also assume that there are at least $m$ distinct points $\{y_{i,j}: i=1,2,...,m\}$ of
$\{z_1,z_2,...,z_N\}$ in each of $\af^j(G)$ (for some $1\le J\le L$), $j=0,1,...,n-1$) such that
\begin{eqnarray}\label{eppm2}
{1\over{n}}\ge {m\over{k}}>{1\over{n}}-{\ep\over{4n}}.
\end{eqnarray}
Furthermore, we may assume $m>L$ and $y_{0,i}\in U_i$ $i=1,2,...,l.$
Put $\Psi(f)=\sum_{i=1}^N f(z_i)p_i$ for $f\in C(X).$
With sufficiently small $\dt_1$ and sufficiently large ${\cal G},$ by
\ref{LMM}, we may also assume
that
\begin{eqnarray}\label{MM1}
\mu_{tr\circ \Psi}(S)\le \mu_{tr\circ \Psi\circ (\af^{-j})^*}(S_{\eta_2/2})
\andeqn
\mu_{tr\circ \Psi\circ(\af^{-j})^*}(S)\le \mu_{tr\circ \Psi}(S_{\eta_2/2})
\end{eqnarray}
for any closed subset $S\subset X,$
where
$(\af^{-j})^*(f)=f\circ \af^{-j},$ $j=1,2,...,n-1$ and
where $S_{\eta_2/2}=\{x\in X: \di(x, S)<\eta_2/2\}.$

Thus, by the choice of $\eta_2,$
 for any $y_{s(i),j},$ $i=1,2,...,M$ with $1\le M\le m,$
there exist $\xi_1',\xi_2',...,\xi_M'\in
\{x\in X: \di(x,\{y_{1,0},y_{2,0},...,y_{m,0}\})<\eta_1/2\}$
such that
$$
\di(y_{s(i),j}
,\af^j(\xi_i'))<\eta_2/2,\,\,\,i=1,2,...,M.
$$
Then, by the choice of $\eta_1,$ there are
$\xi_1,\xi_2,...,\xi_M\in \{y_{1,0},y_{2,0},...,y_{m,0}\}$
such that
$$
\di(y_{s(i),j}
,\af^j(\xi_i))<\eta/2,\,\,\,i=1,2,...,M.
$$
Similarly, for any $\xi_1',\xi_2',...,\xi_M'\in \{y_{1,0},y_{2,0},...,y_{m,0}\},$
there exist $y'_{s(i),j},$ $i=1,2,...,M,$ such that
$$
\di(\af^j(\xi_i'),y'_{s(i),j})<\eta/2\,\,\,\,i=1,2,...,M.
$$
It follows from the ``marriage lemma'' (\cite{HV}) (see also 2.1 of \cite{Su}) that
there is a permutation $\sigma_j: \{1,2,...,m\}\to \{1,2,...,m\}$ such that
$$
\di(y_{i,j}, \af^j(y_{\sigma_j(i),0}))<\eta,
$$
$j=1,2,...,n-1.$
By the choice of $\eta$ and
by replacing $\ep/8$ by $\ep/4$ in (\ref{eppm1}),
we may assume
$y_{i,j}=\af^j(y_{i,0})$
and $y_{i,0}=x_i$ for $1\le i\le l.$
Let $y_i=y_{1,i},$ $i=1,2,...,m.$
Put
$K=mn.$

Thus, from above, with sufficiently large ${\cal G}$ and sufficiently small
$\dt,$ we may also assume that,
\begin{eqnarray}\label{eppm3}
\|\sum_{i=1}^Nf(z_i)p_i-[\sum_{j=0}^{n-1}\sum_{i=1}^m
f(\alpha^j(y_i))p_{i,j}+\sum_{i=K+1}^Nf(z_i)p_i]\|<\ep/2
\end{eqnarray}
for all $f\in {\cal F}.$ Then (\ref{ePms}) follows from (\ref{eppm1}) and (\ref{eppm3}). Moreover, by
(\ref{eppm2}) and (4),
$$
{K\over{k}}={nm\over{k}}>n({1\over{n}}-{\ep\over{4n}})=1-\ep/4.
$$
as desired.




\end{proof}

\begin{Prop}\label{wtre2}
Let $A$ and $B$ be two unital separable \CA s with
$TR(A)=TR(B)=0.$ Suppose that  $\lambda: Aff(T(A))\to Aff(T(B))$ is a
unital order affine isomorphism. Then, there are finite dimensional \CA s $F_n,$ a sequence of unital
\morp s $\phi_n: B\to F_n$ and a sequence of unital \morp s $\psi_n: A\to F_n$ satisfying the
following:

{\rm (1)}
$$
\lim_{n\to\infty}\|\phi_n(a)\phi_n(b)-\phi_n(ab)\|=0
$$
for all $a, b\in A$ and
$$
\lim_{n\to\infty}\|\psi_n(x)\psi_n(y)-\psi_n(xy)\|=0
$$
for all $x,y\in B;$

 {\rm (2)}  there is an affine continuous map $\Delta_n:
T(B)\to T(F_n)$ such that, for each $b\in B,$
\begin{eqnarray}\label{etr1}
|\Delta_n(\tau)(\phi_n(b))-\tau(b)|\to 0
\end{eqnarray}
uniformly on $T(B)$ and

{\rm (3)} for each $a\in A,$
\begin{eqnarray}\label{etr2}
|\lambda(\hat{a})(\tau)-\Delta_n(\tau)\circ\psi_n(a)|\to 0
\end{eqnarray}
uniformly
 on $T(B).$

\end{Prop}

\begin{proof}
Let $\ep>0,$ ${\cal F}\subset A$ and ${\cal G}\subset B$ be two finite subsets. To simplify notation,
without loss of generality, we may assume that ${\cal F}$ and ${\cal G}$ are in the unit balls of $A$
and $B,$ respectively.

Since $TR(A)=0,$ by \cite{Lntr0}, for any $\dt>0,$ there exist a projection $p\in A$ and a finite
dimensional \SCA\, $C$ of $A$ with $p=1_C$ such that

{\rm (i)}\,\, $\|pa-ap\|<\dt/8$ for all $a\in {\cal F};$

{\rm (ii)} \,$\di(pap, C)<\dt/8$ for all $a\in {\cal F}$ and

{\rm (iii)} \, $t(1-q)<\dt/4$ for all $t\in T(A).$

We choose $\dt<\min\{\ep/4,1\}.$ Moreover, by 2.3.5 of \cite{Lnb}, there exists a \morp\, ${\tilde
\psi}': pAp\to C$ such that ${\tilde \psi}(c)=c$ if $c\in C.$ Define ${\tilde \psi}(a)={\tilde
\psi}'(pap)$ for all $a\in A.$

Write $C=\bigoplus_{i=1}^kM_{R(i)}.$ Denote by $e_i$ a minimal rank one projection in $M_{R(i)},$
$i=1,2,...,k.$ Since $TR(B)=0,$ $\rho_B(K_0(B))$ is dense in $Aff(T(B)).$ So there exists a projection
$p_i\in B$ such that
\begin{eqnarray}\label{ePwtre2}
r(\widehat{e_i})(\tau)-\dt/8<\tau(p_i)<r(\widehat{e_i})(\tau)
\end{eqnarray}
for all $\tau\in T(B),$ $i=1,2,...,k.$ Note
$$
\sum_{i=1}^kR(i)[p_i]<[1_B]
$$
in $K_0(B).$ Thus (since $TR(B)=0$) we obtain a \SCA\, $B_0\subset B$ for which there exists an
isomorphism $\psi_1: C\to B_0$ so that $\psi_1(e_i)=p_{i,1},$ $i=1,2,...,k.$

Choose ${\cal G}_1$ which contains ${\cal G}$ and $\psi_1\circ {\tilde \psi}({\cal F})$ as well as a
set of generators of $B_0.$ For any $\dt_1>0,$ there is a projection $q\in B$ and a finite dimensional
\SCA\, $F$ of $B$ with $q=1_F$ such that

{\rm (1)} $\|qb-bq\|<\dt_1/8$ for all $b\in {\cal G}_1;$

{\rm (2)} $\di(qbq, F)<\dt_1/8$ for all $b\in {\cal G}_1$ and

{\rm (3)}$\tau(1-q)<\dt_1/4$ for all $\tau\in T(B).$

We may assume that  $\dt_1<\min\{\ep/4,1\}.$ By 2.3.5 of \cite{Lnb},
we may assume that  there exists a \morp\,
$\phi':qBq\to F$ such that $\phi(b)=b$ if $b\in F.$
Define $\phi: B\to F$ by $\phi(b)=\phi'(qbq)$ for all $b\in B.$ Then $\phi$ is ${\cal
G}_1$-$\dt_1/4$-multiplicative \morp.

Furthermore, by \ref{pertM}, we may assume that there exists a \hm\, $h: B_0\to F$ so that
$$
\|h-\phi|_{B_0}\|<\ep/8.
$$
For each $\tau\in T(B)$ define $\Delta(\tau)={1\over{\tau(q)}}\tau|_F.$ Since, for any $b\in B,$
$$
\tau((1-q)bq)=0=\tau(qb(1-q)),
$$
we have
\begin{eqnarray}\label{ePwtre+1}
|\tau(b)-\tau(qbq)|<\dt_1/4.
\end{eqnarray}
for all $\tau\in T(B).$ With $\dt_1<1,$ for any $f\in F,$
\begin{eqnarray}\label{ePwtre+2}
|\tau(f)-\Delta(\tau)(f)|<(1-{1\over{1-\dt_1/4}})|\tau(f)|<(\dt_1/3)|\tau(f)|
\end{eqnarray}
for all $\tau\in T(B).$ By (2) above, (\ref{ePwtre+1}) and (\ref{ePwtre+2}), we estimate that
\begin{eqnarray}\label{ePwtre}
|\tau(b)-\Delta(\tau)(\phi(b))|<\dt_1/4+\dt_1/8+(\dt_1/3)(1+\dt_1/8)+\dt_1/8<\ep/2
\end{eqnarray}
for all $b\in {\cal G}_1.$

Define $\psi(a)=h\circ ({\tilde \psi}(a)).$ Note that
$\psi$ is from $A$ to $F\subset B$ and it is ${\cal F}$-$\ep$-multiplicative.
We also compute that
$$
|\lambda(\hat{a})(\tau)-\Delta(\tau)(\psi(a))|<\ep
$$
for all $a\in {\cal F}.$

\end{proof}

\section{Uniform approximate conjugacy  in measure}

\begin{Def}\label{Dmeasure}
{\rm Let $X$ be a compact metric space and let $\alpha: X\to X$ be a minimal homeomorphism. Define $F(T_\af)$ to
be  the set of those measures $\nu$ such that
$$
\int f d\nu=\sum_{i=1}^k a_if(x_i),
$$
for all $f\in C(X),$ where $x_1,x_2,...,x_k$ are points in $X,$
$0\le a_i\le 1$ and $\sum_{i=1}^ka_i=1.$

Fix a finite set of points $x_1,x_2,...,x_k\in X$ and $k$ many positive affine continuous functions
$a_1,a_2,...,a_k\in Aff(T(A_\af))$ with $\sum_{i=1}^ka_i=1.$ One can define an affine continuous map
$\Delta: T_\af\to F(T_\af)$ as follows.
\begin{eqnarray}\label{eDel}
\int f d\Delta(\mu)=\sum_{i=1}^k a_i(\tau_{\mu})f(x_i)
\end{eqnarray}
for all $f\in C(X).$
 To simplify notation, we also use $\Delta$
for the induced affine continuous map from $T(A_\af)$ to
$F(T_\af).$

}
\end{Def}

\begin{Def}\label{Dapp1}
{\rm Let $X$ be a compact metric space and  $\af, \bt: X\to X$ be two minimal homeomorphisms. We say
that $\af$ and $\bt$ are {\it approximately conjugate uniformly in measure} if there are
a sequence of open subsets $\{O_n\}$
with  each $O_n$ being $1/n$-dense in $X,$ and a sequence of Borel isomorphisms $\{\gm_n\}$ on $X,$

 {\rm (1)} for each
$\sigma>0,$
\begin{eqnarray}\label{eDap1}
\mu(\{x\in X: \di(\gm_n^{-1}\af\gm_n(x),\bt(x))\ge \sigma\})\to 0
\end{eqnarray}
\begin{eqnarray}\label{eDap1+}
\mu(\{x\in X: \di(\af\gm_n(x),\gm_n\bt(x))\ge \sigma\})\to 0,
\end{eqnarray}
 and
\begin{eqnarray}\label{eDap2}
\nu(\{x\in X: \di(\gm_n\bt\gm_n^{-1}(x), \af(x))\ge \sigma\})\to 0
\end{eqnarray}
\begin{eqnarray}\label{eDap2+}
\nu(\{x\in X: \di(\bt\gm_n^{-1}(x), \gm_n^{-1}\af(x))\ge \sigma\})\to 0
\end{eqnarray}
uniformly on $T_\bt$ and $T_\af,$ respectively,

{\rm (2)} $\gm_n(O_n)$ is a ${1\over{n}}$-dense open subset,  $\gm_n$ is continuous on $O_n$  and
$\gm_n^{-1}$ is  continuous on
$\gm_n(O_n);$

{\rm (3)} there exists an affine continuous map $\Delta_n: T_\bt\to F(T_\bt)$ such that $\int f\circ
\gm_n d\Delta_n(\mu)$ converges uniformly on $T_\bt$ for all $f\in C(X)$ which defines an affine
homeomorphism $r: T_\bt\to T_\af$ and
\begin{eqnarray}\label{eDap3}
|\int f d\mu-\int fd\Delta_n(\mu)|\to 0
\end{eqnarray}
uniformly on $T_\bt$ for all $f\in C(X),$ and there exists an affine continuous map ${\tilde
\Delta}_n: T_\af\to F(T_\af)$ such that $\int f\circ \gm_n^{-1} d{\tilde \Delta}_n(\nu)$ converges
uniformly on $T_\af$ for all $f\in C(X)$ which defines the affine homeomorphism  $r^{-1}:T_\af \to
T_\bt,$ and
\begin{eqnarray}\label{eDap3+}
|\int f d\mu-\int fd{\tilde \Delta}_n(\mu)|\to 0
\end{eqnarray}
uniformly on $T_\af$ for all $f\in C(X).$
 }
\end{Def}

\begin{Remark}\label{RR}
{\rm In general, one should not expect that $\{\gm_n\}$ converges in any suitable sense. Nevertheless,
it is important that $\{\gm_n\}$ carries some consistent information. Note that a homeomorphism does
not preserve measures. Given a sequence of homeomorphisms $\{\gm_n\}$ from $X$ onto $X.$ Even if each
$\gm_n$ does not map positive measure sets to sets with zero measure, it could still happen that, for
example,
$\mu(\gm_n(E))\to 0$ for some Borel set $E$ with $\mu(E)>0.$  Therefore
one should regard (3) as a crucial  part of the definition.

Moreover, given an affine homeomorphism $r: T_\af\to T_\bt,$ Theorem \ref{MT1} provides a sequence of maps
$\{\gm_n\}$ which induces the map $r$ in the sense of (3) in \ref{Dapp1} and
$\gm_n^{-1}\af\gm_n$ converges to $\bt$ and $\gm_n\bt\gm_n^{-1}$
converges to $\af$ uniformly on $T_\bt$ and $T_\af,$ respectively.

 }
\end{Remark}

For convenience, we would like to list two known facts below.

\begin{Lem}\label{lhg}
Let $X$ be a compact metric space and $\af: X\to X$ be a minimal homeomorphism. Then, for any $x, y\in
X$ and any two open neighborhoods
$N(x)$ and $N(y)$ of $x$ and $y,$
there exist a neighborhood $O(x)\subset N(x),$ an open subset
$O\subset N(y)$ and  a homeomorphism $\af'$ from $ O(x)$ onto $O.$
\end{Lem}

\begin{proof}
This follows from the minimality immediately. In fact, for any $\ep>0,$ there exists $n\ge 1,$ such
that
$$
\di(\af^n(x),y)<\ep/2.
$$
Since $\af^n$ is continuous, there exists $\dt>0$ such that
$$
\af^n(\{\xi\in X: \di(x, \xi)<\dt\})\subset \{\xi\in X: \di(y,\xi)<\ep\}.
$$
This means that the homeomorphism $\af^n$ which maps $\{x\in X:\di(x, \xi)<\dt\}$ into the
neighborhood $\{\xi\in X:\di(y,\xi)<\ep\}.$

\end{proof}

\begin{Lem}\label{Borel}
Two second countable locally compact Hausdorff spaces are Borel equivalent if they have the same
cardinality ($\le 2^{\aleph_0}$) .
\end{Lem}

See 4.6.13 of \cite{Ped} for a proof of \ref{Borel}.

\begin{thm}\label{MT1}
Let $X$ be a finite dimensional compact metric space with infinitely  points and let $\af, \bt: X\to
X$ be two minimal homeomorphisms.  Suppose that $\rho(K_0(A_\af))$ is dense in
$Aff(T_\af))$ and $\rho(K_0(A_\bt))$ is dense in $K_0(A_\bt).$
Then the following are equivalent:

{\rm (1)} There is a unital order affine isomorphism $r:Aff(T(A_{\af}))\to Aff(T(A_{\bt}));$

{\rm (2)} $\af$ and $\bt$ are approximately conjugate uniformly in measure;

\end{thm}

\begin{proof}
It suffices to prove  ``{\rm (1)} $\Rightarrow$ {\rm (2)}".

Fix $\ep>0$ and a finite subset ${\cal F}\subset C(X).$ Fix $\eta_0>0$ such that
$$
|f(x)-f(x')|<\ep/8
$$
if $\di(x,x')<\eta_0.$

Choose an integer $n>0$ such that $1/n<\ep/8.$
Choose $\eta_1>0$ such that
\begin{eqnarray}\label{eMP0}
\di(\af^j(x), \af^j(y))<\eta_0/2\andeqn \di(\bt^j(x),\bt^j(y))<\eta_0/2
\end{eqnarray}
if $\di(x,y)<\eta_1,$ $j=1,2,...,n-1.$

Let $\eta=\min\{\ep/4,\eta_1/4,\eta_0/4\}.$


By \ref{Rok1}, one obtains
an open subset $G$  that satisfies the following:
%

{\rm (i)} $G$ contains $\cup_{i=1}^l\{x\in X: \di(x,x_i)<d\}$
for some $d>0,$ where  $\{x_1,x_2,...,x_t\}$ is $\eta/6$-dense;

{\rm (ii)} $\af^j(G)$ are pairwise disjoint for $0\le j\le n-1;$

{\rm (iii)} $\mu(X\setminus \cup_{j=0}^{n-1} \af^j(G))<\ep/8$ for all $\mu\in T_\af$ and

{\rm (iv)} $\mu(\partial(G))=0$ for all $\mu\in T_\af.$




Similarly, let $\Omega$ be an open subset  that  satisfies the following


{\rm (i')} $\Omega$ contains at least one open balls of $\xi_i,$ where $\{\xi_1,\xi_2,...,\xi_t\}$
is $\eta/2$-dense in $X;$


{\rm (ii')} $\bt^j(\Omega)$ are pairwise disjoint for $0\le j\le n-1;$



{\rm (iii')} $\mu(X\setminus \cup_{j=0}^{n-1}\bt^j(\Omega))>1-\ep/8$ for all $\mu\in T_\bt$ and

{\rm (iv')} $\mu(\partial(\Om))=0$ for all $\mu\in T_\af.$

Note that  we can use the same number $t$ for the number of points in
$\{x_1,x_2,...,x_t\}$ and in $\{\xi_1,\xi_2,...,\xi_t\}.$
When we apply \ref{Rok1} to obtain $\Om,$ we use the $\eta/6$-dense set $\{x_1,x_2,...,x_t\}$ to obtain
the $\eta/2$-dense set $\{\xi_1,\xi_2,...,\xi_t\}.$

Suppose that $O(x_i)$ are open balls of $x_i$ so that $O(x_i)\subset G$ and
$O(\xi_i)$ are open balls of $\xi_i$ so that $O(\xi_i)\subset \Omega.$
Since $(X,\af)$ has mean dimension zero, let $\{O_1, O_2,...,O_L\}$ be a finite
set of
pairwise disjoint open subsets
of $X$ such that
each $O_i$ has diameter less than $\eta_1/2,$ $X=\cup_{i=1}^L\overline{O_i}$ and
$\mu(\overline{O_i}\setminus O_i)=0$ for all $\mu\in T_\af.$
We may assume that $O(x_i)\subset O_{i'}\cap G$ for some $i',$ by choosing smaller
open ball of $x_i$ if necessary. Further, by considering a suitable open ball of $x_i$
with universal null boundary,
we may simply assume that $O_i=O(x_i),$ $i=1,2,...,t$ and $L>t.$

Let $\{U_1,U_2,...,U_{L_1}\}$ be a finite set of pairwise disjoint open subsets
of $X$ such that each $U_i$ has diameter less than $\eta_1/2,$
$X=\cup_{i=1}^{L_1}\overline{O_i}$ and
$\nu(\overline{U_i}\setminus U_i)=0$ for all $\nu\in T_\bt.$
We may also assume that $O(\xi_i)=U_{i},$ $i=1,2,...,t$ and $t<L_1.$

Let $p\in A_\af$ and $q\in A_\bt$ be specially selected projections as required by \ref{Permutation}
with
\begin{eqnarray}\label{eMPtrp}
\tau(1-p)<\ep/16 \andeqn \theta(1-q)<\ep/16
\end{eqnarray}
for all $\tau\in T(A_\af)$ and $\theta\in T(A_\bt)$ for
$\ep/4,$ ${\cal F},$ $\eta,$ $n$ and $G$ above and
$\ep/4,$ ${\cal F},$ $\eta,$ $n$ and $\Omega$ above.

Let ${\cal G}_1\subset A_\bt$ be a finite subset (in place of ${\cal G}$) and $\dt>0$ as required by
\ref{Permutation} for the above $\ep/4,$ ${\cal F},$
$n,$ $\eta$ and $\Omega.$ Let ${\cal G}_2\subset A_\af$ be a finite subset
and $\dt_1>0$ as required by  \ref{Permutation} for the above $\ep/4,$ ${\cal F},$
$n$ $\eta$ and $G.$

It follows from \ref{wtre2} (and (\ref{eMPtrp})), with sufficiently large ${\cal G}_1$ and
sufficiently small
$\dt,$ there is a finite dimensional \CA\, $B_0,$
a unital  ${\cal G}_1$-$\dt$-multiplicative \morp\, $\phi: A_{\beta} \to B_0,$ a ${\cal
G}_2$-$\dt$-multiplicative \morp\, $\psi: A_{\alpha}\to  B_0$ and an affine continuous map $\Delta_0:
T(A_\bt) \to T(B_0)$ such that

{\rm (1)}
\begin{eqnarray}\label{eMP1}
|\Delta_0(\tau)\circ \phi(pj_{\beta}(f)p)-\tau\circ j_{\beta}(f)|<\ep/8
\end{eqnarray}
for all $\tau\in T(A_{\beta})$ and $f\in {\cal F}$ and

{\rm (2)}
\begin{eqnarray}\label{eMP2}
|r(\widehat{j_{\alpha}(f)})(\tau)-\Delta_0(\tau)\circ\psi(pj_{\alpha}(f)p)|<\ep/8
\end{eqnarray}
for all $\tau\in T(A_{\beta})$ and $f\in {\cal F}.$

Write $B_0=\oplus_{s=1}^{k_0}M_{R(s)}$ and $\pi_s: B_0\to M_{R(s)}$ the canonical projection map.
By applying \ref{Permutation},  for
each $s,$ there are integers $K(s)=m_sn$ and $K'(s)=m_s'n$ with
$m_s=\sum_{i=1}^Lm_s(i)$ and $m_s'=\sum_{i'=1}^{L_1}m_s'(i'),$
there are points $y_{i,l}(s)\in O_i\cap G,$ $l=1,2,...,m_s(i)$ and
$i=1,2,...,L,$ $Y_{i',l'}(s)\in U_i\cap \Om,$ $l'=1,2,...,m_s'(i')$
and $i'=1,2,...,L_1$ such that
\begin{eqnarray}\label{eMP4}
\|\sum_{i,l,j}f(\af^j(y_{i,l}(s)))p_{s,i,l,j}+\sum_{i=K(s)+1}^{N(s)}f(z_i)p_{s,i}-\pi_s\circ\psi\circ
(pj_\af(f)p)\|<\ep/4
\end{eqnarray}
for all $f\in {\cal F}$ and
\begin{eqnarray}\label{eMP5}
\|\sum_{i',l',j}f(\beta^j(Y_{i,l}(s)))q_{s,i',l',j}+
\sum_{i'=K'(s)+1}^{N'(s)}f(z'_i)q_{s,i'}-\pi_s\circ\phi\circ (qj_\bt(f)q)\|<\ep/4
\end{eqnarray}
for all $f\in {\cal F},$ where
$$
\{p_{s,i,l,j}: i,l,j\}\cup \{p_{s,i}: i>N(s)\} \andeqn \{q_{s,i',l',j}: i',l',j
\}\cup\{q_{s,i'}:i'>N'(s)\}
$$
are sets of mutually orthogonal rank one projections in $M_{R(s)}$ and $z_i, z_{i'}\in X.$ In
addition, by \ref{Permutation}, we may assume that $y_{i,1}(1)=x_i$  and
$Y_{i,1}(1)=\xi_i$ $i=1,2,...,t.$

Furthermore,
\begin{eqnarray}\label{eMPINT}
{R(s)-K(s)\over{R(s)}}<\ep/4\andeqn {R(s)-K'(s)\over{R(s)}}<\ep/4,
\end{eqnarray}
$s=1,2,...,k_0.$ Without loss of generality, since $X$ has no
isolated points, we may assume that $\{y_{i,l}(s): i,l,s\}$ and
$\{Y_{i',l'})(s): i',l',s\}$ are two sets of distinct points. If
$m_s'>m_s,$ we will move $m_s'-m_s$ many points of $Y_{i',l'}(s)$
to the set $\{z_i': i'\}.$  If, on the other hand, $m_s>m_s',$ we will move $m_s-m_s'$ many points to
$\{z_i: i\}.$ So,  either way, we may assume that $m_s=m_s'$ and $K(s)=K'(s).$ Note that we still have
${R(s)-K'(s)\over{R(s)}}<\ep/4.$

By replacing $\phi$ by ${\rm ad}\, u\circ \phi,$ for a suitable unitary in $B_0,$ we may assume that
$$
\{p_{s,i,l,j}: 1\le i\le L,1\le l\le m_s(i),0\le j\le n-1\}=\{q_{s,i',l',j}\}.
$$

Since now we assume that $m_s=m'_s,$ we define, for each $s,$
${\tilde \gm}(Y_{i',l'}(s))$ to be an one to one bijection
between $\{Y_{i'l'}(s): i',l',s\}$ and $\{y_{i,l}(s): i,l,s\}.$
We may also assume that ${\tilde \gm}(Y_{i,1}(1))=y_{i,1}(1),$
$i=1,2,...,t.$

To construct the desired map $\gm,$ we divide $O_i\cap G$ into $\sum_{s=1}^{k_0}m_s(i)$ pairwise
disjoint sets $B_{s,i,l}$ as follows: $B_{s,i,1}$ is an open subset which contains an open
neighborhood of $y_{i,1}(s)$ for
$1\le i\le t$ and every  other $B_{s,i,l}$ is the closures of an open
neighborhood of $y_{i,l}(s)$ ($1\le l\le m_s$). We then divide
 $U_i\cap\Om$ into
$\sum_{s=1}^{k_0}m_s'(i')$ pairwise disjoint  subsets  $C_{s,i',l'}$ as follows:
$C_{s,i',1}$ is an open subset which contains
an open neighborhood of $Y_{i',1}(s)$ for $1\le i'\le t.$ Every other
$C_{s,i',l'}$ is  the closure of an open neighborhood of
$Y_{i',l'}(s)$ ($1\le l'\le m_s'$).

Note that, since every points in $X$ is condensed, $B_{s,i,l}$ and $C_{s,i',l'}$ are second countable
locally compact Hausdorff spaces with the cardinality $2^{\aleph_0}.$ By \ref{Borel}, they are all
Borel equivalent.

Define a Borel equivalence $\gm: X\to X$ as follows:

By \ref{lhg}, there is  an open neighborhood $Z(i,1,s)$ of $Y_{i,1}(s)$ in $C_{s,i,1}$ (for $1\le i\le
t$) and  a open subset ${\tilde Z}(i,1,s)$ of  $B_{s,i,1}$ which are homeomorphic. In particular, the
closure of a smaller open neighborhood of $Y_{i,1}(s)$ is homeomorphic to a closure of an open subset
of ${\tilde Z}(i,1,s).$
 Thus, by taking sufficiently small  such neighborhood and
by applying \ref{Borel},  one obtains a Borel equivalence $\gm$ from $C_{s,i',1}$ onto
$B_{s,i,1}$ which maps  a non-empty
neighborhood $Z(i,1,s)$ of $Y_{i,1}(s)$ to an open subset of a neighborhood   of
$y_{i,1}(s)$ homeomorphically for $1\le i\le t.$

We define $\gm$ on $\bt^j(C_{s,i',l'})$ to be $\af^j\circ \gm\circ \bt^{-j},$
$j=1,2,..,n-2.$

Since
 $X\setminus
\cup_{j=0}^{n-2}\af^j(G)$ (which is a compact subset of $X$ containing
$\af^{n-1}(G)$) and $X\setminus \cup_{j=0}^{n-2}\af^j(\Om_{j})$ (which is a compact subset
of $X$ which contains
$\af^{n-1}(\Om)$) are
Borel equivalent,
 we obtain a Borel equivalence $\gm$ of $X$ which is a bi-continuous on
$O=\cup_{i',s} Z(i',1,s).$ Note that $\gm$ maps
$\bigcup_{j=0}^{n-2}\bt^j(\Om)$ onto
$\bigcup_{j=0}^{n-2}\af^j(G).$
We also have
$\gm(Z(i,1,s))\subset {\tilde Z}(i,1,s).$
 Since $\cup_{i=1}^L O_i$ and
$\cup_{i'=1}^{L_1}U_i$
have  the diameter  less than $\eta/2,$ from the construction, we see that $O$ and
$\gm(O)$ are $\eta$ -dense in $X.$

Moreover, on each $\bt^j(C_{s,i'l'})$ with $0\le j\le n-2,$
\begin{eqnarray}\label{eMP9}
\di(\gm^{-1}\af\gm(x),\bt(x))<\eta\andeqn \di(\af\gm(x),\gm\bt(x))<\eta
\end{eqnarray}
We also have, on each $\af^j(B_{i,l,s})$ with $0\le j\le n-2,$
\begin{eqnarray}\label{eMP10}
\di(\gm\bt \gm^{-1}(x),\af(x))<\eta\andeqn \di(\bt\gm^{-1}(x),\gm^{-1}\af(x))<\eta
\end{eqnarray}

Since
\begin{eqnarray}\label{eMP11}
\nu(\bt^{n-1}(\Om))<1/n<\ep/8\andeqn \mu(\af^{n-1}(G))<1/n<\ep/8
\end{eqnarray}
for all $\bt$-invariant probability measure $\nu$ and $\af$-invariant probability measure
$\mu,$
we conclude that
\begin{eqnarray}\label{eMP12}
\nu(\{x\in X: \di(\gm^{-1}\alpha\gm(x),\bt(x))>\eta\})<\ep/4
\andeqn\\
\mu(\{x\in X: \di(\gm\bt\gm^{-1}(x), \af(x))>\eta\})<\ep/4
\end{eqnarray}
for all $\bt$-invariant probability measure $\nu$ and $\af$-invariant probability measure
$\mu.$

To complete the proof, it remains to check (3) of \ref{Dapp1}. For this end, we note, by (\ref{eMP4}),
(\ref{eMP5}) and (\ref{eMPINT}), that
\begin{eqnarray}\label{eMPtrab}
|\sum_{s,i,l,j}f(\af^j(y_{i,l}(s)))\Delta_0(\tau)(p_{s,i,l,j})-\Delta_0(\tau)(\psi\circ
(pj_\af(f)p))|<\ep/2
\end{eqnarray}
 and
\begin{eqnarray}\label{eMPtrab2}
|\sum_{s,i',l',j'}f(\beta^j(Y_{i,l}(s)))\Delta_0(\tau)(q_{s,i',l',j}) -\Delta_0(\tau)(\phi\circ
(qj_\bt(f)q))|<\ep/2
\end{eqnarray}
for all $f\in {\cal F}$ and $\tau\in T(A_\bt).$ Note also, for each $s,$
$\Delta_0(\tau)(p_{s,i,l,j})=\Delta_0(\tau)(q_{s,i',l',j})={c_{\tau}\over{R(s)}}$ for all $i,i',l,l',j$
and for some non-negative constant $c_{\tau}.$

We also estimate that, for each $s,$
\begin{eqnarray}\label{eMP13}
|\sum_{i',l', 0\le j\le n-2}f\circ \gm(\beta^j(Y_{i',l'}(s))){c_{\tau}\over{R(s)}}- \sum_{i,l,0\le
j\le n-2} f(\af^j(y_{i,l}(s))){c_{\tau}\over{R(s)}}|<\ep/8
\end{eqnarray}
and
\begin{eqnarray}\label{eMP13+}
|\sum_{i,l,o\le j\le n-2}f\circ \gm^{-1}(\af^j(y_{i,l}(s))){c_{\tau}\over{R(s)}}- \sum_{i,l,0\le j\le
n-2} f(\bt^j(Y_{i',l'}(s))){c_{\tau}\over{R(s)}}|<\ep/8
\end{eqnarray}
for all $f\in {\cal F}$ and $\tau\in T_\bt.$

Define $\Delta: T_\bt \to F(T_\bt)$ by
$$
\int fd(\Delta(\mu))=\sum_{s=1}^{k_0}\sum_{i',l',0\le j\le n-2}f(\bt^j(Y_{i'.l'}(s))){c_{\tau}\over{R(s)}}
$$
for all $\mu\in T_\bt$ (where $\mu=\mu_\tau$) and all $f\in C(X).$
Note that $\int fd(r_{\natural}(\mu))=r(\widehat{j_{\af}(f)})(\tau)$ ($\mu=\mu_{\tau}$).
Combining  (\ref{eMPtrp}), (\ref{eMP1}), (\ref{eMP2}), (\ref{eMPINT}), (\ref{eMPtrab})  and
(\ref{eMP13}), we have
\begin{eqnarray}\label{eMP14}
|\int f d\mu-\int f d(\Delta(\mu)|<\ep,\\
|\int f\circ \gm d(\Delta(\mu))-\int f d(r_{\natural}(\mu))|<\ep
\end{eqnarray}
for all $\bt$-invariant probability measure $\mu$ and all $f\in {\cal F}.$

Define ${\tilde \Delta}: T_\af\to F(T_\af)$ by
${\tilde \Delta}(\nu)=\Delta(r_{\natural}^{-1}(\nu))$ for $\nu\in T_\af.$
Then, we have, by (\ref{eMPtrp}),  (\ref{eMP1}), (\ref{eMP2}), (\ref{eMPINT}), (\ref{eMPtrab2})  and
(\ref{eMP13+}),
\begin{eqnarray}\label{eMP15}
|\int f d{\tilde \Delta}(\nu)-\int fd\nu|<\ep\\
|\int f\circ \gm^{-1} d{\tilde \Delta}(\nu)-\int f d r_{\natural}^{-1}(\nu)|<\ep
\end{eqnarray}
for all $f\in {\cal F}.$

\end{proof}

\section{Concluding remarks}



({\bf 1})\,\,\, Let $X$ be a compact metric space and $T$ be a convex subset of probability Borel
measures. Suppose that
$\Gamma_n, \Gamma: X\to X$ are Borel maps and $\Gamma_n\to
\Lambda_n$ in measure uniformly on $T.$ Then a uniform E{\tiny
$\Gamma$}o$\rho$o{\tiny B} theorem holds. Put
\begin{eqnarray}
S_{m,k}=\{x\in X: \di(\Gamma_m, \Gamma(x))\ge 1/k\},
\end{eqnarray}
$k=1,2,...,$ and $m=1,2,....$ Let $\dt>0.$
For each $k>0,$ there exists an integer $n(k)$ such that
\begin{eqnarray}\label{eCR1}
\mu(S_{n(k),k})<{\dt\over{2^k}}
\end{eqnarray}
for all $\mu\in T,$  if $n\ge n(k).$ Put
\begin{eqnarray}
E=\bigcap_{k=1}^{\infty}\bigcap_{m=n(k)}^{\infty}\{x\in X: \di(\Gamma_m(x),\Gamma(x))<1/k\}.
\end{eqnarray}
Then $\Gamma_n$ converges to $\Gamma$ uniformly on $E.$ Furthermore,
\begin{eqnarray}\label{eCR2}
\mu(X\setminus E)\le \mu(\bigcup_{k=1}^{\infty}S_{n(k),k})\le \sum_{k=1}^{\infty} \mu(S_{n(k),k})<\dt
\end{eqnarray}
for all $\mu\in T.$ Thus, in Theorem \ref{MT1}, for any $\dt>0,$ there exists a Borel subset $E\subset
X$ with $\mu(X\setminus E)<\dt$ for all $\mu\in T_\bt$ such that $\gm_n^{-1}\af\gm_n$ converges to
$\bt$ uniformly on $E.$  Moreover, there exists a Borel subset
$E'\subset X$ with $\mu(X\setminus E')<\dt$ such that
$\gm_n\bt\gm_n^{-1}$ converges to $\af$ uniformly on $E'.$
A similar measure theoretical argument, by taking a subsequence of $\{\gm_n\},$ shows that there exist
Borel measurable subsets $F_\af, F_\bt \subset X$ such that
$\gm_n^{-1}\af\gm_n$ converges to $\bt$  on $F_\bt$ and
$\gm_n\bt\gm_n^{-1}$ converges to $\af$ on $F_\af$ and
$X\setminus F_\bt $  and $X\setminus F_\bt$ are universally null, i.e., $\mu(X\setminus F_\bt)=0$
for all $\mu\in T_\bt$ and $\nu(X\setminus F_\af)=0$ for all $\nu\in T_\af.$


\vspace{0.2in}

\noindent ({\bf 2})\,\,\,
Suppose that $X$ is the Cantor set and suppose that $\af,\bt: X\to X$ are
two minimal homeomorphisms. Then in Theorem \ref{Rok1} $G$ can be chosen to be
clopen.
Since a non-empty
clopen subset of the Cantor set can be divided into $m$ non-empty clopen subsets for any integer
$m>0,$ in the proof of \ref{MT1}, $B_{i,l,s}$ and
$C_{i',l',s}$ can be chosen to be also non-empty clopen subsets of
$X.$ They all are homeomorphic. It is then easy to see that the
map $\gm$ in the proof can be made to be homeomorphism. In other words, we have the following
corollary:

\begin{Cor}\label{CCantor}
Let $X$ be the Cantor set and let $\af, \bt: X\to X$ be minimal homeomorphisms. Then $\af$ and $\bt$
are approximately conjugate in measure if and only if there is an affine homeomorphism $r: T_\af$ to
$T_\bt.$ Moreover, when $\af$ and $\bt$ are approximately conjugate uniformly in measure, the
conjugating maps
$\gm_n$ can be chosen to be homeomorphisms.
\end{Cor}




\end{document}